\documentclass[12pt,a4paper]{article}

\usepackage{hyperref}
\usepackage[english]{babel} 
\usepackage{amsthm}
\usepackage{amsmath}
\usepackage{amsfonts}
\usepackage{amssymb}
\usepackage{enumerate}
\usepackage{amscd}
\usepackage{graphicx}
\usepackage{makeidx}
\usepackage[all]{xy}

\makeindex

\newcounter{Def}[section]
\renewcommand{\theDef}{\arabic{section}.\arabic{Def}}
 
\theoremstyle{plain}
\newtheorem{Proposition}[Def]{Proposition}

\newtheorem{Examples}[Def]{Examples}
\newtheorem{Example}[Def]{Example}
\newtheorem{Remarks}[Def]{Remarks}
\newtheorem{Remark}[Def]{Remark}
\newtheorem{Lemma}[Def]{Lemma}
\newtheorem{Theorem}[Def] {Theorem}
\newtheorem{Kor}[Def] {Corollary}

\newtheorem{corollary}[Def]{Corollary}

\newenvironment{Proof}[1][\quad]
{\mbox{}\\\noindent \textbf{Proof:{\hspace{0.7cm}
#1}\\}}{\hfill$\Box$\\}

\newenvironment{Definition}{ \refstepcounter{Def}\mbox{}\\ 
\noindent\sl\textbf{Definition
\arabic{section}.\arabic{Def}} }
{\vspace{0.3cm}}

\def\BiCat         {\ensuremath{\mathcal{B}\!i\mathcal{C}\!at}}

\def\C          {\mathbb C}

\DeclareMathOperator{\holim}{holim}
\def\Z          {\mathbb Z}

\def\calo        {\mathcal O}

\def\calg         {{\mathcal G}}
\def\calj         { {\mathcal J}}
\def\Ce            {\text{\v{C}}\,}
\def\be                 {\begin{equation}}
\def\ee                 {\end{equation}}
\def\Desc          {\mathcal{D}\!esc}

\def\gto           {\rightrightarrows}

\def\LieGrp              {\mathrm{LieGrpd}}
\def\Kan                 {\mathfrak{C}\mathrm{an}}

\def\oti           {\,{\otimes}\,}

\def\id               {{\rm id}}
\def\ii          {\mathrm{i}}
\def\I                 {\text{I}}

\def\To            {\,{\to}\,}

\def\rmd         {\mathrm{d}}

\def\id                 {\mathrm{id}}

\def\X                  {\mathfrak{X}}

\newcommand{\iso}{\stackrel{\sim}{\longrightarrow}}

\def\c                  {\mathcal{C}}  
\def\d                  {\mathcal{D}} 
\def\Hom{\textnormal{Hom}}
\def\Man           {\ensuremath{\mathcal M\mbox{\sl an}}}
\def\uto           {\twoheadrightarrow}
\def\BiCat         {\ensuremath{\mathcal B\mbox{\sl i} \mathcal C \mbox{\sl at}}}
\def\CTo           {\twoheadrightarrow}
\def\calg          {{\cal G}}
\def\Grbc          {\Grb^{\!\nabla}\!}
\def\Grb           {\mathcal{G}\hspace{-.8pt}\mbox{\sl rb}}
\def\Grbtriv      {\Grb\mbox{\sl triv}}
\def\Grbctriv      {\Grb\mbox{\sl triv}^{\!\nabla}\!}
\def\JGrbc          {\JGrb^{\!\nabla}\!}
\def\JOGrbc         {\JGrb_{\text{or }}^{\!\nabla}\!}
\def\JGrb           {\mathcal{JG}\hspace{-.8pt}\mbox{\sl rb}}

\def\JGrbctriv      {\JGrb\mbox{\sl triv}^{\!\nabla}\!}
\def\hol           {\mathrm{Hol}}
\def\Ue            {\mathrm U(1)}
\def\curv          {\mathrm{curv}}
\def\Bun           {\mathcal{B}{\!un}}
\def\Bunc          {\Bun\!^\nabla\!}
\def\JBun           {\mathcal{JB}{\!un}}
\def\JBunc          {\JBun\!^\nabla\!}
\def\Buntriv       {\mathcal{B}{\!untriv}}

\def\vect          {\textnormal{Vect}_\mathbb{C}}
\def\zvecttriv      {\mathcal{V}{\!ect_2 triv}}
\def \zvect         {\mathcal{V}{\!ect_2}}

\begin{document}

\numberwithin{equation}{section}


\thispagestyle{empty}
\begin{flushright}
   {\sf ZMP-HH/10-12}\\
   {\sf Hamburger$\;$Beitr\"age$\;$zur$\;$Mathematik$\;$Nr.$\;$373}\\[2mm]
   April 2010
\end{flushright}
\vskip 2.0em
\begin{center}\Large
EQUIVARIANCE\ IN\ HIGHER\ GEOMETRY
\end{center}\vskip 1.4em
\begin{center}
  Thomas Nikolaus and Christoph Schweigert\,\footnote{\scriptsize 
  ~Email addresses: \\
  $~$\hspace*{2.4em}Thomas.Nikolaus@uni-hamburg.de, 
  Christoph.Schweigert@uni-hamburg.de}
\end{center}

\vskip 3mm

\begin{center}\it
  Fachbereich Mathematik, \ Universit\"at Hamburg\\
  Bereich Algebra und Zahlentheorie\\
  Bundesstra\ss e 55, \ D\,--\,20\,146\, Hamburg
\end{center}
\vskip 2.5em
\begin{abstract} \noindent
We study (pre-)sheaves in bicategories on geometric categories:
smooth manifolds, manifolds with a Lie group action and Lie groupoids.
We present three main results: we describe equivariant descent,
we generalize the plus construction to our setting and
show that the plus construction yields a 2-stackification for 2-prestacks.
Finally we show that, for a 2-stack, the pullback functor along a 
Morita-equivalence of Lie groupoids is an equivalence of bicategories.

Our results have direct applications to gerbes and 2-vector bundles.
For instance, they allow to construct equivariant gerbes from local data 
and can be used to simplify the description of the local data. We
illustrate the usefulness of our results in a systematic discussion of
holonomies for unoriented surfaces.

\end{abstract}


\setcounter{footnote}{0} \def\thefootnote{\arabic{footnote}} 

\tableofcontents

\section{Introduction}\label{sec:introduction}

In a typical geometric situation, one
selects a category of geometric spaces, e.g.\ smooth
manifolds, and then considers
for every geometric space $M$ a category $\mathfrak{X}(M)$
of geometric objects on $M$, e.g.\ complex line bundles or
principal $G$-bundles, with $G$ a Lie group.
The categories for different geometric spaces
are related by pullback functors: they
form a presheaf in categories.

In this paper, the category of geometric spaces we consider
is the category $\LieGrp$ of Lie groupoids. This category has crucial advantages:
it contains \v Cech groupoids and thus provides a convenient
setting to discuss local data. Moreover, it contains
action groupoids and thus allows us to deal with equivariant
geometric objects as well.

We show that any presheaf $\X$ on manifolds can be naturally
extended to a presheaf on Lie groupoids.
We also generalize the structure we associate to
a geometric space $M$ by considering a bicategory $\X(M)$.
This choice is motivated by the fact that
bundle gerbes and bundle gerbes with connection on a
given manifold have the structure of a bicategory
\cite{stev2000, waldorf2007more}.
Hence we will work with a 
presheaf in bicategories on the geometric category $\LieGrp$
of Lie groupoids. Our theory extends the theory for 
(pre-)sheaves in categories on smooth manifolds presented
in \cite{metzler2003topological,heinloth2005notes}.

A hallmark of any geometric theory is a procedure to obtain
global objects from locally defined objects by a gluing procedure.
To this end, one considers open covers which are, in the category
of smooth manifolds, just a special class $\tau_{open}$ of morphisms.
More generally, we endow the category of manifolds  with
a Grothendieck topology, although we will not directly use
this language to keep this article at a more elementary
level. The two prime examples for choices of $\tau$
for the category of smooth manifolds are $\tau_{open}$, i.e.\ open
covers, and $\tau_{sub}$, i.e.\ surjective submersions.

Having fixed a choice for $\tau$, we get a notion of
$\tau$-essential surjectivity of Lie functors and of
$\tau$-weak equivalence of Lie groupoids $\Gamma$ and $\Lambda$.
($\tau_{sub}$-weak equivalent Lie groupoids are also called {\em Morita equivalent};
some authors also call a $\tau_{sub}$-weak equivalence a Morita equivalence.)
Imposing different gluing conditions on the presheaf
$\X$ on $\LieGrp$ for morphisms in $\tau$,
we get the notion of a $\tau$-2-prestack on $\LieGrp$ and
of a $\tau$-2-stack on $\LieGrp$, respectively. To simplify the notation,
we refer to a 2-prestack as a prestack
and to a 2-stack as a stack.

These basic definitions are the subject of section \ref{sec:sheafLie}. At the 
end of this section, we
can state our first main theorem \ref{2.13}:

\setcounter{section}{2}
\setcounter{Def}{15}

\begin{Theorem} \mbox{} \label{1.1} \\
Suppose, $\Gamma$ and $\Lambda$ are Lie groupoids and
$\Gamma\to\Lambda$ is a $\tau$-weak equivalence of
Lie groupoids.
\begin{enumerate}
\item
Let $\mathfrak X$ be a $\tau$-prestack on $\LieGrp$.
Then the functor
$$\ \mathfrak{X}(\Lambda) \to \mathfrak{X}(\Gamma) $$
given by pullback is fully faithful, i.e.\ an equivalence on the Hom categories.

\item
Let $\mathfrak X$ be a $\tau$-stack on $\LieGrp$.  Then the functor
$$\ \mathfrak{X}(\Lambda) \to \mathfrak{X}(\Gamma) $$
given by pullback is an equivalence of bicategories.
\end{enumerate}
\end{Theorem}

This theorem, or more precisely its first assertion, is a central
ingredient for our second main result which we explain in section 
\ref{sec:plus}. In analogy to the sheafification of
a presheaf, we associate to any prestack $\X$ a presheaf in bicategories 
$\X^+$ where
the objects of the bicategory $\X^+(M)$ consist of a cover $Y\to M$ and an
object in the descent bicategory $\Desc_\X(Y\to M)$. We call this construction
the plus construction. We then state theorem \ref{stackifizierung}:

\setcounter{section}{3}
\setcounter{Def}{2}

\begin{Theorem}\mbox{} \label{1.2} \\
Let $\X$ be a prestack on $\Man$. Then
the presheaf in bicategories $\X^+$ on $\Man$ obtained by the plus construction 
is a stack. Furthermore the canonical
embedding $\X(M) \to \X^+(M)$ is fully faithful for each manifold $M$.
\end{Theorem}

The plus construction is a powerful tool to construct geometric objects. 
In section \ref{sec:applplus}, we show this in the example of bundle gerbes 
with connection: we introduce a bicategory $\Grbctriv$ of trivial bundle 
gerbes with connection whose objects are given by 2-forms. 
A brief check reveals that the
plus construction yields bundle gerbes,
$$ \Grbc = \big( \Grbctriv\big)^+ \,\,\, . $$
Theorem \ref{stackifizierung} then immediately
implies that bundle gerbes form a stack.

Bundle gerbes give rise to a notion of surface holonomy.
We then apply the reasoning leading to the definition of bundle gerbes
to the definition of surface holonomy
for unoriented surfaces and find the notion of a Jandl gerbe. In appendix 
\ref{sapp:sholunor}, we also compare this notion to the notion of a Jandl
structure on a gerbe that has been introduced earlier \cite{ssw}. 
Based on the notion of Jandl gerbe, we introduce in appendix 
\ref{sapp:sholunor} the notion of an 
orientifold background on a Lie groupoid $\Lambda$.  Theorem
\ref{2.13} allows us to define a surface holonomy for any
Hilsum-Skandalis morphism \cite[definition 62]{metzler2003topological}
from the unoriented worldsheet $\Sigma$ to
$\Lambda$.  

It should be stressed that our results apply to general higher geometric
objects, in particular to non-abelian gerbes and 2-vector bundles. To
illustrate this point, subsection \ref{ssec:kvvect} contains a short discussion
of 2-vector bundles. In all cases, theorem \ref{1.2} immediately ensures
that these higher geometric objects form a stack over the category of 
manifolds (and even of Lie groupoids).

Together, these results provide us with tools to construct concrete geometric
objects: theorem \ref{1.2} allows us to glue together geometric
objects like e.g.\ gerbes from locally defined geometric object.
Applications frequently require not only gerbes, but equivariant gerbes.
Here, it pays off that our approach is set off for Lie groupoids rather than
for manifolds only, since the latter combine equivariance and local data
on the same footing. In particular, we are able to formulate in this
framework theorem \ref{HauptsatzEquivDesc} on {\em equivariant} descent. One application
of this theorem is to obtain equivariant gerbes from locally defined
equivariant gerbes.

Theorem \ref{1.1} and theorem \ref{HauptsatzEquivDesc} can then
be combined with standard results on the action of
Lie groups or Lie groupoids \cite{Duistermaat-Kolk,Weinstein}
to obtain a simplified description of the local situation
in terms of stabilizer groups. This strategy provides,
in particular, an
elegant understanding of equivariant higher categorical geometric objects,
see e.g.\
\cite{nikolaus-thesis} for the construction of
gerbes on compact Lie groups \cite{Meinrenken, Gawedzki-Reis}
that are equivariant under the adjoint action.

We have collected the proofs of the theorems in the second part
of this paper in sections \ref{sec:proof1} -- \ref{sec:plusproof}. In 
an appendix, we discuss applications to surface holonomies and systematically
introduce a notion of holonomy for unoriented surfaces.

\bigskip

\noindent{\bf Acknowledgements.}
We thank Till Barmeier, Urs Schreiber and Konrad Waldorf for helpful
discussions and Lukas Buhn{\'e} and Konrad Waldorf for comments on the draft.
TN and CS are partially supported by the Collaborative Research Centre 676 
``Particles, Strings and the Early Universe - the Structure of Matter 
and Space-Time''.

\setcounter{section}{1}
\section{Sheaves on Lie groupoids}\label{sec:sheafLie}

\subsection{Lie groupoids}\label{ssec:LieGrp}
We start our discussion with an introduction to Lie groupoids. Groupoids
are categories in which all morphisms are isomorphisms. A small
groupoid, more specifically, consists of a set $\Gamma_0$ of objects
and a set $\Gamma_1$ of morphisms, together with maps
$s,t: \Gamma_1 \rightarrow \Gamma_0, \iota: \Gamma_0 \rightarrow \Gamma_1$ that associate to a morphism $f\in\Gamma_1$ its source
$s(f)\in\Gamma_0$ and its target $t(f)\in \Gamma_0$ and to
an object $m\in\Gamma_0$ the identity $\id_m\in\Gamma_1$.
Finally, there is an involution $in:\Gamma_1\to\Gamma_1$
that obeys the axioms of an inverse.
Concatenation is a map 
$\circ: \Gamma_1 \times_{\Gamma_0} \Gamma_1 \rightarrow \Gamma_1$
where it should be appreciated that in the category of sets
the pullback $\Gamma_1 \times_{\Gamma_0} \Gamma_1
=\{(f_1,f_2)\in\Gamma_1\times\Gamma_1 | t(f_1)=s(f_2)\}$
exists. It is straightforward to translate the usual axioms
of a category into commuting diagrams.

A Lie groupoid is groupoid object in the category of smooth manifolds:

\begin{Definition}\mbox{} \\
A groupoid in the category $\Man$ or a  {\em Lie-groupoid} 
consists of two smooth manifolds $\Gamma_0$ and $\Gamma_1$ 
together with the following collection of smooth maps:
\begin{itemize}
\item Source and target maps 
$s,t: \Gamma_1 \rightarrow \Gamma_0$.
\end{itemize}
To be able to define compositions, we need the existence of the
pullback $\Gamma_1 \times_{\Gamma_0} \Gamma_1$. To ensure 
its existence, we require $s$ and $t$ to be surjective submersions. 

The other structural maps are:
\begin{itemize}
\item 
A composition map $\circ: \Gamma_1 \times_{\Gamma_0} \Gamma_1 \rightarrow 
\Gamma_1$ 
\item A neutral map $\iota: \Gamma_0 \rightarrow \Gamma_1$ providing identities
\item A map $in: \Gamma_1\to \Gamma_1$ giving inverses
\end{itemize}
such that the usual diagrams commute.
\end{Definition}

\begin{Examples} \mbox{}\label{2.2} \\[-1.8em]
\begin{enumerate}
\item For any manifold, we have the trivial Lie groupoid
$M \rightrightarrows M$ in which all structure maps are
identities. We use this to embed $\Man$ into $\LieGrp$.

\item Given any Lie group $G$, we consider the Lie groupoid $BG$
with structure maps $G \rightrightarrows$pt with pt the smooth zero-dimensional
manifold consisting of a single point. The neutral map
$pt \rightarrow G$ is given by the map to the neutral element and composition
$G \times G \rightarrow G$ is group multiplication.
Hence Lie groupoids are also a generalization of Lie groups.

\item More generally, if a Lie group $G$ is acting smoothly
on a smooth manifold $M$, the \emph{action groupoid}
$M // G$ has $\Gamma_0:=M$ as objects and the manifold
$\Gamma_1:=G \times M$ as morphisms. The source map $s$
is projection to $M$, the target map $t$ is given by the
action $t(g,m) := g\cdot m $. The neutral map is the injection
$m \mapsto (1,m)$ and composition is given by the group
product, $(g,m) \circ (h,n) := (gh, n)$. Action Lie groupoids
frequently are the appropriate generalizations of
quotient spaces.

\item For any covering $(U_i)_{i\in I}$ of a manifold $M$
by open sets $U_i\subset M$, we consider the disjoint 
union $Y:=\sqcup_{i\in I} U_i$ with the natural local homeomorphism
$\pi: Y \twoheadrightarrow M$. Consider the two natural
projections $Y \times_M Y \rightrightarrows Y$ with the
composition map
$(Y \times_M Y) \times_Y (Y \times_M Y) \cong Y^{[3]} 
\to Y^{[2]}$ given by omission of the second element.
The neutral map is the diagonal map
$Y \to Y \times_M Y$. This defines a groupoid $\Ce(Y)$,
the \v Cech-groupoid.

\end{enumerate}
\end{Examples}

The last two examples show that Lie groupoids provide a
convenient framework to unify ``local data'' and 
equivariant objects.

We next need to introduce morphisms of Lie groupoids.

\begin{Definition} \mbox{} \\
A morphism of Lie groupoids or \emph{Lie functor}
$F: (\Gamma_1 \rightrightarrows \Gamma_0) 
\to (\Omega_1 \rightrightarrows \Omega_0)$ 
consists of smooth maps $F_0: \Gamma_0 \to \Omega_0$ and $F_1: \Gamma_1 \to \Omega_1$ that are required to commute with the structure
maps. For example, for the source map $s$, we have the
commuting diagram
$$\xymatrix{
\Gamma_1 \ar[r]^{F_1}\ar[d]_{s} & \Omega_1 \ar[d]^s\\
\Gamma_0 \ar[r]^{F_0}           & \Omega_0
}$$

\end{Definition}

\begin{Examples} \mbox{}\label{2.4} \\[-1.8em]
\begin{enumerate}
\item Given two smooth manifolds $M,N$, every Lie functor
$F: (M \rightrightarrows M) \to (N \rightrightarrows N)$
is given by a smooth map $f: M \to N$  with $F_0 = F_1 =f$. 
Hence  $M \mapsto (M\rightrightarrows M)$ is a fully faithful
embedding and we identify the manifold $M$ with the
Lie groupoid $M\rightrightarrows M$.

\item 
Given two Lie groups $G$\ and $H$, the Lie functors
$F:BG\to BH$ between the corresponding Lie groupoids are
given by smooth group homomorphisms $f: G \to H$.
Thus the functor $G \mapsto BG$ is a fully faithful embedding of
Lie groups into Lie groupoids.

\item For any two action groupoids $M // G$ and $N // G$,
a $G$-equivariant map $f: M \to N$ provides a Lie functor
via $F_0 := f$ and $F_1 := f \times id : M \times G 
\to N \times G$. The previous example with $M=N=pt$ shows that
not all Lie functors between action groupoids are of this 
form.

\item 
Consider a refinement $Z \uto M$ of a covering
$Y \uto M$ together with the refinement map $s: Z \to Y$.
This provides a Lie functor $S:\Ce(Z)\to\Ce(Y)$
of \v Cech groupoids which acts on objects by
$S_0 := s : Z \to Y$ and on morphisms 
$S_1 : Z \times_M Z \to Y \times_M Y$ by $S_1(z_1,z_2) 
:= (s(z_1),s(z_2)) \in Y \times_M Y$. 

\item  As a special case, any covering $Y \uto M$ is a refinement
of the trivial covering $id: M \uto M$ and we obtain a Lie
functor $\Pi^Y: \Ce(Y) \to M$.

\end{enumerate}
\end{Examples}
\subsection{Presheaves in bicategories on Lie groupoids}\label{ssec:prebi}

A presheaf in bicategories $\X$ on the category $\Man$ of manifolds consist of a 
bicategory \cite{benabou1967}
$\X(M)$ for each manifold $M$, a pullback functor $f^*: \X(N) \to \X(M)$
for each smooth map $f: M \to N$ and natural isomorphisms $f^* \circ g^* \cong (g \circ f)^*$
for composable smooth maps $f$ and $g$. Moreover, we need higher coherence isomorphisms satisfying 
the obvious, but lengthy conditions. More precisely, $\X$ is a weak functor
$$ \X: \Man^{op} \to \BiCat.$$
Furthermore we impose the technical condition that $\X$ preserves products, i.e. for a disjoint union $M = \bigsqcup_{i \in I} M_i$ of manifolds indexed by a set $I$ the following equivalence holds:
\begin{equation}\label{disjoint_union}
 \X(M) \cong \prod_{i \in I} \X(M_i) \,\, .
\end{equation}

Our next step is to extend such a presheaf in bicategories on $\Man$ to a 
presheaf 
in bicategories on Lie groupoids. For a Lie groupoid $\Gamma$
finite fiber products $\Gamma_1\times_{\Gamma_0}\cdots
\times_{\Gamma_0}\Gamma_1$ exist in $\Man$ and we
introduce the notation
$\Gamma_2=\Gamma_1\times_{\Gamma_0}\Gamma_1$
and $\Gamma_n$ analogously.

We can then use the nerve construction to associate to
a Lie groupoid a simplicial manifold
\begin{equation*}
\left(\xymatrix{
{\cdots}  
\ar@<1.3ex>[r]^{\partial_0}
\ar@<0.3ex>[r]
\ar@<-0.7ex>[r]
\ar@<-1.7ex>[r]_{\partial_3}
&
{\Gamma_2}
\ar@<0.9ex>[r]^{\partial_0}
\ar@<-0.1ex>[r]
\ar@<-1.1ex>[r]_{\partial_2}
&
{\Gamma_1}
\ar@<0.3ex>[r]^{\partial_0}
\ar@<-0.7ex>[r]_{\partial_1} 
&
{\Gamma_0}
}\right)=:\Gamma_\bullet \,\,\, .
\end{equation*}
We can think of $\Gamma_n$ as $n$-tuples of morphisms
in $\Gamma_1$ that can be concatenated. The map 
$\partial_i: \Gamma_{n} \rightarrow \Gamma_{n-1} $
is given by composition of the $i$-th and $i+1$-th morphism.\ Thus

\begin{eqnarray*} 
\partial_i(f_1,\ldots,f_n) &:=& (f_1,\ldots,f_i \circ f_{i+1},\ldots,f_n) \\
\partial_0(f_1,\ldots,f_n) &:=& (f_2,\ldots,f_n) \\
\partial_n(f_1,\ldots,f_n) &:=& (f_1,\ldots,f_{n-1}) \,\, .
\end{eqnarray*}
In particular, $\partial_1,\partial_0: \Gamma_1 \rightarrow \Gamma_0$ 
are the source and target map of the groupoid. One easily verifies
the simplicial identities
$\partial_i \partial_{j+1} = \partial_j \partial_i$ for 
$i \leq j$. (We suppress the discussion of the degeneracy maps
$\sigma_i: \Gamma_n \rightarrow \Gamma_{n+1}$ which are given
by insertion of an identity morphism at the $i$-the position.)

The nerve construction can also be applied to Lie functors and
provides an embedding of Lie groupoids into simplicial manifolds.
Suppose we are given a Lie functor  
$F: (\Gamma_1 \rightrightarrows \Gamma_0) \to (\Omega_1\rightrightarrows \Omega_0)$.
Consider the nerves $\Gamma_\bullet$ and $\Omega_\bullet$
and define a family $F_\bullet= (F_i)$ of maps, a
{\em simplicial map}
$$F_i: \Gamma_i \to \Omega_i$$
for all $i=0,1,2,\ldots$ with $F_0,F_1$ given by the Lie functor
and maps given for $i >1 $ by
\begin{eqnarray*}
     F_i: \quad\Gamma_1 \times_{\Gamma_0} \ldots \times_{\Gamma_0} \Gamma_1
         &\to& \Omega_1 \times_{\Omega_0} \ldots \times_{\Omega_0} \Omega_1 \\
      (f_1,\ldots,f_n) &\mapsto& \big(F_1(f_1),\ldots F_1(f_n)\big)\,\, .
\end{eqnarray*}
By definition, the maps $F_i$ commute with the maps
$\partial_j$ and $\sigma_k$ that are part of the simplicial object.
We summarize this in the following diagram:
\begin{equation*}
\xymatrix{
{\cdots}  
\ar@<1.3ex>[r]
\ar@<0.3ex>[r]
\ar@<-0.7ex>[r]
\ar@<-1.7ex>[r]
&
{\Gamma_2}
\ar@<0.9ex>[r]
\ar@<-0.1ex>[r]
\ar@<-1.1ex>[r]
\ar[d]_{F_2}
&
{\Gamma_1}
\ar@<0.3ex>[r]
\ar@<-0.7ex>[r]
\ar[d]_{F_1}
&
{\Gamma_0}
\ar[d]_{F_0} 
\\
{\cdots}  
\ar@<1.3ex>[r]
\ar@<0.3ex>[r]
\ar@<-0.7ex>[r]
\ar@<-1.7ex>[r]
&
{\Omega_2}
\ar@<0.9ex>[r]
\ar@<-0.1ex>[r]
\ar@<-1.1ex>[r]
&
{\Omega_1}
\ar@<0.3ex>[r]
\ar@<-0.7ex>[r]
&
{\Omega_0}
}
\end{equation*}

\begin{Definition}\label{defEquivariant} \mbox{} \\
Let $\mathfrak{X}$ be a presheaf in bicategories on $\Man$
and $\Gamma$\ a Lie groupoid or, more generally, a simplicial
manifold. A $\Gamma$-{\em equivariant object} of $\mathfrak X$ consists of 
\begin{itemize}
\item[(O1)]
an object $\mathcal{G}$ of $\mathfrak{X}(\Gamma_0)$;
\item[(O2)]
a 1-isomorphism
  \begin{equation*}
  P\colon~ \partial_0^{*}\mathcal{G} \To \partial_1^{*}\mathcal{G}
  \end{equation*}
in $\X(\Gamma_1)$;
\item[(O3)]
a 2-isomorphism
  \begin{equation*}
  \mu \colon~ \partial_{2}^*P \oti \partial_{0}^*P \,{\Rightarrow}\, \partial_{1}^*P
  \end{equation*}
in $\X(\Gamma_2)$, where we denote the horizontal product by $\otimes$;
\item[(O4)]
a coherence condition
 \begin{equation*}
  \partial_2^*\mu \circ (\id \oti \partial_0^* \mu) = \partial_1^*\mu \circ (\partial_3^*\mu \oti \id)\,
  \end{equation*}
on 2-morphisms in $\X(\Gamma_3)$.
\end{itemize}
\end{Definition}

We next introduce 1-morphisms and 2-morphisms of $\Gamma$-equivariant objects:

\begin{Definition}  \mbox{} \label{2.6} \\[-1.8em]
\begin{enumerate}
\item
A 1-morphism between $\Gamma$-equivariant objects
$(\mathcal{G}, P, \mu)$ and $(\mathcal{G}', P', \mu')$ in
$\X$ consists of the following data on the simplicial manifold
  $$
  \left(\ldots
\xymatrix{\Gamma_4 \ar[r]<.6em>\ar[r]<.2em>\ar[r]<-.2em>\ar[r]<-.6em> &
  \,\Gamma_2 \ar[r]<.4em>\ar[r]<.0em>\ar[r]<-.4em> &
  \,\Gamma_1 \ar_-{\partial_1}[r]<-.2em> \ar^-{\partial_0}[r]<.2em> &
  \,\Gamma_0
  }\right)=\Gamma_\bullet
  $$
\def\leftmargini{3.91em}
\begin{itemize}
\item[(1M1)] 
A 1-morphism $A\colon \mathcal{G} \to \mathcal{G}'$ in $\X(\Gamma_0)$;
\item[(1M2)]
A 2-isomorphism $\alpha\colon P' \oti \partial_0^*A \,{\Rightarrow} \,\partial_1^*A \oti P$ 
in $\X(\Gamma_1)$;
\item[(1M3)]
A commutative diagram
  \begin{equation*}
  (\id \oti \mu) \circ (\partial_2^*\alpha \oti \id) \circ (\id \oti \partial_0^*\alpha)
  = \partial_1^* \alpha \circ (\mu' \oti id)
  \end{equation*}
of 2-morphisms in $\X(\Gamma_2)$.
\end{itemize}

\item
A 2-morphism between two such 1-morphisms
$(A,\alpha)$ and $(A',\alpha')$ consists of
\def\leftmargini{3.85em}\\[-1.5em]
\begin{itemize}
\item[(2M1)]
A 2-morphism $\beta\colon A \, {\Rightarrow} \, A'$ in $\X(\Gamma_0)$;
\item[(2M2)]
a commutative diagram 
  \begin{equation*}
  \alpha' \circ (\id \oti \partial_0^*\beta) = (\partial_1^*\beta \oti \id) \circ \alpha
  \end{equation*}
of 2-morphisms in $\X(\Gamma_1)$.
\end{itemize}
\end{enumerate}
\end{Definition}

We define the composition of morphisms using simplicial identities and
composition in the bicategories $\X(\Gamma_i)$, see e.g.\ \cite{waldorf2007more}.
The relevant definitions are lengthy but
straightforward, and we refrain from giving details.

One can check that in this way, one obtains the structure of a 
bicategory.

\begin{Remarks}\label{2.7} \mbox{}\\[-2em]
\begin{enumerate}
\item Similar descent bicategories have been introduced
in \cite{breen1994classification} and \cite{duskin1989outline}.
For a related discussion of equivariance in presheaves in bicategories,
see also \cite{skoda}. 

\item
We have defined $\Gamma$-equivariant objects  for  a presheaf
$\X$ in bicategories. Any presheaf $\X$ in categories can be
considered as a presheaf in bicategories with trivial
2-morphisms. We thus obtain a definition for $\X(\Gamma)$
for presheaves in categories as well, where the 2-morphisms in
(O3) on $\Gamma_3$ become identities and the condition
(O4) is trivially fulfilled. Similar remarks apply to
morphisms. All 2-morphisms are identities, hence 
$\mathfrak{X}(\Gamma_\bullet)$ can be identified with a
category. This allows us to deal with presheaves in
categories as special cases of our more general
results on presheaves in bicategories and to recover part of the results of \cite{metzler2003topological,heinloth2005notes}.
\end{enumerate}
\end{Remarks}

One can check that the following proposition holds:

\begin{Proposition} \mbox{}\label{2.8} \\
Our construction provides for any Lie groupoid $\Gamma$ a bicategory
$\X(\Gamma)$. The bicategories form a presheaf in bicategories
on the category $\LieGrp$ of Lie groupoids.

\end{Proposition}

To make contact with existing literature, we introduce
for the special case of an action groupoid $N//G$ as in example
\ref{2.2}.3
the alternative notation
$$\X_G(N) := \X(N//G)\text{.}$$

\begin{Remarks} \mbox{} \label{rem:explicit} \\[-1.8em]
\begin{enumerate}
\item 
For the convenience of the reader, we spell out the definition of a
$G$-equivariant object of a presheaf in bicategories $\X$ for the special 
case of a
discrete group $G$. A $G$-equivariant object on a $G$-manifold $N$
consists of
\begin{itemize}
\item An object $\calg\in\X(N)$.
\item For every group element $g\in G$ a morphism
$g^*\calg\stackrel{\varphi_g}\to \calg$.
\item A coherence 2-isomorphism for every pair of group elements $g,h\in G$,
$$\xymatrix{
g^*h^*\calg \ar[rr]^{g^* \varphi_h}\ar[rrd]_{\varphi_{hg}} & & g^* \calg\ar[d]^{\varphi_g} \\
&  &\calg
\ar@{=>} (25,-3)*{}; (20,-8)*{}
}
$$
\item
A coherence condition.
\end{itemize}

\item
We also show how to obtain the usual definition of equivariant bundles
on a $G$-manifold $N$, where $G$ is a Lie group. We denote the action by
$w: N\times G\to N$. An equivariant bundle on $N$ consists of the following 
data: a bundle $\pi: P\to N$ on $N$. The simplicial map
$\partial_0: N\times G \to N$ is projection, $\partial_1=w$ is the action.
Hence $\partial_0^*P= P\times G$ and $\partial_1^*P=w^*P$. The second data
is a morphism 
$P\times G \to w^*P = (N\times G)\times_N P $. A morphism to a fibre product
is a commuting diagram
$$ \xymatrix{
P\times G \ar[r] \ar[d]& P \ar^\pi[d]\\
N\times G \ar^-w[r] & N 
} $$
The left vertical map is bound to be $\pi\times \id_G$. The coherence
condition of the equivariant object tells us that $\tilde w:P\times G\to P$
is in fact a $G$-action that covers the $G$-action on $N$.
\end{enumerate}
\end{Remarks}

\begin{Kor}\mbox {} \\
Let $G$\ be a Lie group. Then the functor
$\X_G$ forms a presheaf in bicategories on the category 
$\Man_G$ of smooth manifolds with $G$ action. 
\end{Kor}

By abuse of notation, we denote the presheaf in bicategories 
on $\LieGrp$ introduced in proposition \ref{2.8}
by $\X$. This is justified by the fact that for a constant
Lie groupoid $M \rightrightarrows M$ one has the equivalence 
$\X(M \rightrightarrows  M) \cong \X(M)$.

We next wish to impose generalizations of the sheaf
conditions on a presheaf. To this end, we have to single out a
collection $\tau$ of 
morphisms in $\Man$. Technically, such a collection 
should form a Grothendieck (pre-)topology. 
This means essentially that the collection $\tau$ of morphisms 
is closed under compositions, pullbacks and contains all identities. 
See \cite{metzler2003topological} for a detailed introduction. 
 For our purposes, two families are important:\ 
\begin{itemize}
\item
The family $\tau_{sub}$ of surjective submersions.
\item
The family $\tau_{open}$ that consists of morphisms
obtained from an open covering $(U_i)_{i\in I}$ of a manifold 
$M$ by taking the local homeomorphism
$\pi:\,\, Y\twoheadrightarrow M$
with $Y:=\sqcup_{i\in I} U_i$.
\end{itemize}

{}From now on, two-headed arrows will be reserved
for morphisms in the relevant topology $\rho$. 
Whenever, in the sequel, no explicit topology is mentioned, we refer
to $\tau_{sub}$ as our standard (pre-)topology.

\begin{Remarks} \mbox{} \\
Let $\rho$ be a topology on $\Man$.
\begin{itemize}
\item For any morphism
$\pi:Y\uto M$ in $\rho$, we can form a \v{C}ech groupoid 
$\Ce(Y)$ as in example \ref{2.2}.4 which we again call
the \v Cech groupoid.

\item 
Given a morphism $\pi:Y\uto M$ of $\rho$, we 
define the descent bicategory  by
$$\Desc_\X(Y \uto M):= \X\big(\Ce(Y)\big) \,\,\, .$$

Recall the Lie functor $\Pi^Y:\,\,\check C(Y)\to M$ for the \v Cech cover
$Y\uto M$ introduced in example \ref{2.4}.5. Applying the presheaf
functor $\X$ to this Lie functor, gives the functor of bicategories
\be \tau_Y:\ \X(M) \to  \X(\check C(Y))=\Desc_\X(Y\uto M) 
\label{eq:2.9}\ee
\end{itemize}
\end{Remarks}

We are now ready for two definitions:

\begin{Definition}\mbox{} \\
Let $\X$ be a presheaf in bicategories on $\Man$ and
$\tau$\ a topology on $\Man$.
\begin{enumerate}
\item 
A presheaf $\X$ is called a $\tau$-\em{prestack}, if for
every covering $Y\uto M$ in $\tau$ the functor
$\tau_Y$ of bicategories in (\ref{eq:2.9}) is fully 
faithful. (A functor of bicategories is called
fully faithful, if  all functors on Hom categories are
equivalences of categories.)

\item 
A presheaf $\X$ is called a $\tau$-\em{stack}, if for
every covering $Y\uto M$ in $\tau$ the functor
$\tau_Y$ of bicategories is an equivalence of bicategories.
\end{enumerate}
\end{Definition}

Generalizing the discussion of \cite[Section 5.4]{momr} for submersions,
we use the topology $\tau$ to single out certain
morphisms of Lie groupoids that we call $\tau$-weak equivalences of
Lie groupoids. To motivate our definition, we discuss
equivalences of small categories $\c,\d$. A functor 
$F:\c\to\d$ is an equivalence, if it is fully faithful
and essentially surjective. The latter condition means
that for any object $d\in\d$, there exists an object
$c\in\c$ and a isomorphism $F(c)\stackrel{f}\to d$ in
$\d$. If the category $\d$ is a groupoid, this amounts to the requirement
that the map from
$$\ \c_0\times_{\d_0}\d_1
=\{ (c,f) | c\in\c_0=Ob(\c),
f\in \d_1=Mor(\d) \text{ with }
F(c)=s(f)\ \} $$
to $\d_0$ induced by the target map is surjective. In the context
of Lie groupoids, we will require this map to be in $\tau$.

\begin{Definition}\label{volltreu}\mbox{} \\[-1.8em]
\begin{enumerate}
\item
A morphism of Lie groupoids $\Gamma\to \Lambda$ is
called \em{fully faithful}, if the diagram
$$\xymatrix{
\Gamma_1\ar^{F_1}[r]\ar_{s\times t}[d]&\Lambda_1 \ar_{s\times t}[d]\\
\Gamma_0\times \Gamma_0 \ar_{F_0\times F_0}[r]&
\Lambda_0\times\Lambda_0
}$$
is a pull back diagram.

\item
A morphism of Lie groupoids $\Gamma\to \Lambda$ is
called $\tau$-{\em essentially surjective}, if the
smooth map
$$\ \Gamma_0\times_{\Lambda_0}\Lambda_1 \to \Lambda_0\ $$
induced by the target map in $\Lambda$ is in $\tau$.

\item A Lie functor is called a $\tau$-weak equivalence of Lie
groupoids, if it is fully faithful and $\tau$-essentially
surjective. If we omit the prefix $\tau$, we always 
refer to $\tau_{sub}$-weak equivalences.
\end{enumerate}
\end{Definition}

\begin{Remark} \mbox{} \\
Despite its name, $\tau$-weak equivalence is not an equivalence relation.
The equivalence relation generated by $\tau_{sub}$-weak equivalences is called
Morita equivalence or, for general $\tau$-weak equivalences $\tau$-Morita
equivalence. Explicitly, two Lie groupoids $\Gamma$ and $\Lambda$ are
Morita equivalent, if there exists a third Lie groupoid $\Omega$ and
$\tau$-weak equivalences $\Gamma\to \Omega$ and $\Lambda\to \Omega$.
\end{Remark}

\begin{Example}\mbox{}\label{bsp:cech} \\
The Lie functor $\Pi^Y:\,\,\check C(Y)\to M$ is a $\tau$-weak equivalence for
all $\tau$-covers.
\end{Example}

The stack axiom just asserts that for all 
$\tau$-weak equivalences of this type, the induced functor 
on bicategories
$ \tau_Y:\ \X(M) \to  \X(\check C(Y))$ is an equivalence of
bicategories. The first theorem of this paper
generalizes this statement to all $\tau$-weak equivalences
of Lie groupoids:

\begin{Theorem} \mbox{}\label{2.13} \\
Suppose, $\Gamma$ and $\Lambda$ are Lie groupoids and
$\Gamma\to\Lambda$ is a $\tau$-weak equivalence of
Lie groupoids. 
\begin{enumerate}
\item
Let $\mathfrak X$ be a $\tau$-prestack on $\LieGrp$. 
Then the functor
$$\ \mathfrak{X}(\Lambda) \to \mathfrak{X}(\Gamma) $$
given by pullback is fully faithful.

\item
Let $\mathfrak X$ be a $\tau$-stack on $\LieGrp$.  Then the functor
$$\ \mathfrak{X}(\Lambda) \to \mathfrak{X}(\Gamma) $$
given by pullback is an equivalence of bicategories.
\end{enumerate}
\end{Theorem}

Roughly speaking,  $\tau$-covers of manifolds can be thought
of as being
dense enough in $\tau$-weak equivalences of Lie groupoids to 
allow an extension of 
the (pre-)stack condition.

We defer the proof of the theorem to section \ref{sec:proof1}
- \ref{sec:proof4}
and first present some applications.

\subsection{Open coverings versus surjective submersions}\label{sec:vs}

We have already introduced two Grothendieck (pre-)topologies
$\tau_{open}$ and $\tau_{sub}$ on the category of smooth
manifolds. Since open covers are special examples of
surjective submersions, any $\tau_{sub}$-(pre)stack is 
obviously a $\tau_{open}$-(pre)stack. From theorem \ref{2.13},
we deduce the converse:

\begin{Proposition}\mbox{} \\
A presheaf in bicategories on $\LieGrp$ is a 
$\tau_{open}$-(pre)stack if and only if it is a
$\tau_{sub}$-(pre)stack.
\end{Proposition}

The proposition implies in particular that it is enough to check
the stack condition on open covers.

\begin{Proof}
It remains to be shown that any $\tau_{open}$-stack
$\X$ is also a $\tau_{sub}$-stack. We fix a surjective
submersion $\pi:Y\uto M$ and obtain a functor
$$ \tau_Y:\,\,\, \X(M)\ \to \Desc_\X(Y\uto M)=\X(\Ce(Y)) \,\, . $$
For the surjective submersion $\pi$, we can find local sections
$$\ s_i:\ U_i \to Y\ $$
for an open cover $(U_i)_{i\in I}$ of $M$. We glue together
these sections to a map $s$ on the disjoint union of the open
subsets. Then the diagram
$$\xymatrix{
\sqcup_{i\in I} U_i\ar^s[r]\ar[dr]&\ Y\ar^\pi[d] \\
&M
} $$
commutes. Here the unlabeled arrow is the inclusion of open subsets.
This diagram induces a commuting diagram of Lie groupoids
$$\xymatrix{
\Ce(\sqcup_{i\in I} U_i)\ar^s[r]\ar[dr]&\ \Ce(Y)\ar^\pi[d] \\
&M
} $$
in which $s$ is an $\tau_{open}$-weak equivalence of Lie groupoids.
Since $\X$ is a $\tau_{open}$-stack, the application of
$\X$ yields a diagram that commutes up
to a 2-cell,
$$\xymatrix{
\Desc_\X(\sqcup_{i\in I} U_i)&\ar_<<<<{s^*}[l] \Desc_\X(Y) \\
&\X(M)\ar_{\pi^*}[u]\ar[ul]
} $$
We wish to show that the vertical arrow is an equivalence
of bicategories. The lower left arrow is an equivalence of bicategories,
since $\X$ is assumed to be a $\tau_{open}$-stack.
Since $s$ is a $\tau_{open}$-weak equivalence
of Lie groupoids, theorem \ref{2.13} implies that
$s^*$ is an equivalence of bicategories and the assertion 
follows.
\end{Proof}

Since presheaves in categories are particular examples,
an immediate corollary is:

\begin{Kor}\mbox{}\label{3.2} \\
A presheaf in categories on $\LieGrp$ is a 
$\tau_{open}$-(pre)stack if and only if it is a
$\tau_{sub}$-(pre)stack.
\end{Kor}

After one further decategorification, we also obtain

\begin{Kor}\mbox{} \\
A presheaf  on $\LieGrp$ is a 
$\tau_{open}$-separated presheaf if and only if it is a
$\tau_{sub}$-separated presheaf. \\
A presheaf  on $\LieGrp$ is a 
$\tau_{open}$-sheaf if and only if it is a
$\tau_{sub}$-sheaf.
\end{Kor}

Let us discuss an application of this result:
$U(1)$ principal bundles form a stack on $\Man$
with respect to the open topology $\tau_{open}$,
see e.g. \cite{ecm}. As a consequence of corollary
\ref{3.2}, $U(1)$ bundles also form a stack with respect
to surjective submersions. Hence we can glue bundles
also with respect of surjective submersions. In this way, we
recover the following well-known

\begin{Proposition} \mbox{} \\
Consider a free action groupoid $M//G$ such that the quotient 
space $M/G$ has a natural structure of a smooth manifold and the
canonical projection is a submersion. (This is, e.g., the case if
the action of $G$ on $M$ is proper and discontinuous.)
Then the category of smooth $U(1)$-bundles on
$M/G$ is equivalent to the category of $G$-equivariant
$U(1)$-bundles on $M$.
\end{Proposition}

\begin{Proof}
Since the action is free, the canonical projection 
$\pi:\ M\to M/G$ is a submersion that induces a 
$\tau_{sub}$-weak equivalence of Lie groupoids. We have seen that
$U(1)$-bundles form a $\tau_{sub}$-stack, and hence
by theorem \ref{2.13} the canonical projection $\pi$ 
induces an equivalence of categories.
\end{Proof}

We have formulated this result for the special case of 
$U(1)$ bundles. Obviously, the same argument applies to
any stack on $\Man$, and we obtain similar equivalences
of categories for $G$-equivariant principal bundles, 
and associated bundles for any structure group.

\section{The plus construction}\label{sec:plus}
 
In this section we describe a general procedure for 2-stackification.
More precisely,
we show how to obtain a 2-stack $\X^+$ on $\Man$ starting 
from 2-prestack $\X$ on $\Man$. In analogy to the case of sheaves,
we call this construction the plus construction. The idea is to 
complement the bicategories $\X(M)$ by adding objects in descent 
bicategories. The main result is then that the 2-presheaf in bicategories
obtained in this way is closed under descent.  

We first describe the bicategory $\X^+(M)$ for a manifold $M$.

\begin{Definition}\label{xplus}\\
An object of $\X^+(M)$ consists of a covering $Y \twoheadrightarrow M$ 
and an object $G$ in the descent bicategory $\Desc_\X(Y)$.
\end{Definition}

In order to define 1-morphisms and 2-morphisms between objects with possibly 
different coverings $\pi: \, Y \CTo M$ and $\pi': \, Y' \CTo M$, we
pull all the data back to a common refinement of these coverings and compare them 
there. We call a covering $\zeta: Z \CTo M$ a \emph{common refinement} of $\pi$ and $\pi'$ 
iff there exist coverings $s\colon Z \CTo Y$ and $s'\colon Z \CTo Y'$ such that the diagram 
  \be
  \xymatrix{
  Y \ar@{->>}[dr]_-{\pi}  & Z \ar@{->>}[l]_-{s} \ar@{->>}[r]^{s'}\ar@{->>}[d]_{\zeta}& Y'\ar@{->>}[dl]^-{\!\!\pi'} \\
  &M&
  }
  \label{eq:3.1}
  \ee
commutes. An example of such a common refinement is the fibre product  
$Z \,{:=}\, Y {\times_M}\, Y' \CTo M$, with the maps $Z \CTo Y$ and 
$Z \CTo Y'$ given by the projections. We call this the \emph{canonical common refinement}.
The maps $s$ and $s'$ of a common refinement $Z\uto M$ induce Lie functors on the Lie groupoids
$$\xymatrix{  
{\Ce(Y)} & \Ce(Z) \ar[l] \ar[r] & \Ce(Y')\text{.}
}$$
Hence we have \emph{refinement functors} $s^*$ and $s'^*$:
$$\xymatrix{  
\Desc_\X(Y) \ar^-{s^*}[r]& \Desc_\X(Z)   
& \Desc_\X(Y') \ar_-{(s')^*}[l]\text{.}
}$$
For an object $G$ in $\Desc_\X(Y)$ we denote the refinement $s^*(G)$ by $G_Z$.

\begin{Definition}\label{morxplus}
\begin{itemize}
\item A \emph{1-morphism} between objects $\calg = (Y,G)$ and 
$\calg' = (Y',G')$ of $\X^+(M)$ consists of a common refinement $Z \CTo M$ of the 
coverings $Y \CTo M$ and $Y'\CTo M$ 
and a 1-morphism $A: G_Z \to G'_{Z}$ of the two refinements in $\Desc_\X(Z)$. 
\item A 2-morphism between 1-morphisms $\mathfrak{m} = (Z,A)$ and 
$\mathfrak{m}' = (Z',A')$ consists of a common refinement $W\CTo M$ of 
the coverings $Z \CTo M$ and $Z' \CTo M$ 
and a 2-morphism $\beta: \mathfrak{m}_W \Rightarrow  \mathfrak{m}_W'$ of the refined morphisms 
in $\Desc_\X(W)$. In addition two such 2-morphisms
$(W,\beta)$ and $(W',\beta')$ must be identified iff there exists a further 
common refinement $V \CTo M$ of $W \CTo M$ and $W' \CTo M$ such that the refined 2-morphisms agree on $V$. 
\end{itemize}
\end{Definition}

Now that we have defined objects, morphisms and 2-morphisms in $\X^+(M)$ it remains to define compositions and identities. We will just indicate how this is done. For example let $\calg = (Y,G)$, $\calg' = (Y',G')$ and $\calg'' = (Y'',G'')$ be objects and $\mathfrak{m} = (Z,A) : \calg \to \calg'$ and $\mathfrak{m}' = (Z',A'): \calg' \to \calg''$ be morphisms. The covers can then be arranged to the diagram
$$
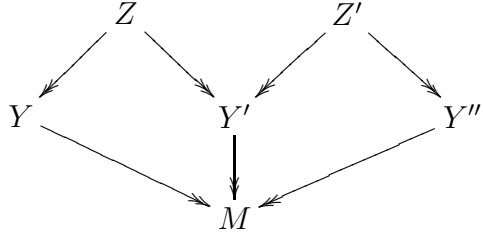
\xymatrix{
  & Z \ar@{->>}[ld]\ar@{->>}[rd]&    & Z'\ar@{->>}[ld]\ar@{->>}[rd] &     \\
Y\ar@{->>}[rrd] &   & Y'\ar@{->>}[d] &    & Y''\ar@{->>}[lld] \\
  &   & M  &    &
}
$$
Now let $Z'' := Z \times_{Y'} Z'$ be the pullback of the upper diagram. 
This exists in $\Man$ and is evidently a common refinement of $Y$ and $Y''$. 
The composition $\mathfrak{m}' \circ \mathfrak{m}$ is then
defined to be the tuple $(Z'', A'_{Z''} \circ A_{Z''})$ where $A'_{Z''} \circ A_{Z''}$ denotes the composition of the refined morphisms 
in $\Desc_\X(Z'')$. \\

Finally one can check that this defines the structure of a bicategory $\X^+(M)$. See 
\cite{waldorf2007more} for a very detailed treatment of a related bicategory. In order to turn the bicategories $\X^+(M)$ into a stack we have to define the pullback functors
$$ f^*: \X^+(N) \to \X^+(M)$$
for all smooth maps $f: M \to N$. This is done in the obvious way using the pullback of covers and the pullback functors of the prestack $\X$.

\begin{Theorem}\label{stackifizierung}\mbox{} \\ 
If $\X$ is a prestack, then $\X^+$ is a stack. Furthermore the canonical 
embedding $\X(M) \to \X^+(M)$ is fully faithful for each $M$. 
\end{Theorem}

We relegate the proof of this theorem to section \ref{sec:plusproof}.

\begin{Remark} \mbox{} \\[-1.8em]
\begin{enumerate}
\item
If we choose the covers in definition \ref{xplus} and \ref{morxplus} to be in the topology $\tau_{open}$ we obtain a slightly different stack $\X^+_{open}$. Argument similar to the ones used in section \ref{sec:vs}
show that  $\X^+(M) \cong \X^+_{open}(M)$ for each smooth manifold $M$.

\item
As in remark \ref{2.7}.2, one can specialize to presheaves in
categories and obtains the stackification process for 1-prestacks.
\end{enumerate}
\end{Remark}

\section{Applications of the plus construction}\label{sec:applplus}
\subsection{Bundle gerbes}\label{gerbes}
We next present several applications of the plus construction.
The input for this construction is a presheaf in bicategories on $\Man$.
In the same way a monoid is the simplest example of a category
(with one object), any monoidal category gives a bicategory with
a single object. An example for a bicategory can thus be 
obtained from the monoidal category
of principal $A$-bundles, where $A$ is any abelian
Lie group. This way, we get a presheaf $\Grbtriv_A$ of of trivial $A$-gerbes. 
Since bundles can be glued together, the homomorphism categories
are closed under descent. The presheaf $\Grbtriv_A$ is thus
a prestack. The plus construction yields the stack 
$$ \Grb_A := \big(\Grbtriv_A)^+ $$
of gerbes (without connection). Our general result
implies that gerbes form a sheaf on $\Man$. Together with
theorem \ref{2.13} and theorem \ref{HauptsatzEquivDesc} of 
this paper, this provides a local construction of gerbes 
and the definition of equivariant gerbes.

Let us next construct gerbes with connection; for simplicity, we restrict
to the abelian group $A=U(1)$ and suppress the index $A$. The guiding principle for
our construction is the requirement that gerbes should lead to a notion
of surface holonomy (for a review, see e.g. \cite{ecm}
and appendix \ref{sapp:sholor}). Hence we consider
a bicategory whose objects are two-forms. These two-forms are intended to
be integrated; hence it is not natural to require them to be equal, but
only to require them to be equal up to a total derivative or, in physical
terms, up to a gauge transformation. We are thus 
lead to consider for every manifold $M$ the following bicategory

\begin{itemize}
\item
Objects are in bijection with 2-forms $\omega \in \Omega^2(M)$ and
denoted by $\mathcal{I}_{\omega}$. 

\item
A 1-morphism $\mathcal{I}_\omega \to \mathcal{I}_{\omega'}$
is a 1-form $\lambda$ such that $\rmd\lambda=\omega'-\omega$.

\item
A 2-morphism $\lambda \to \lambda'$ is a $U(1)$-valued function 
$g$ on $M$ such that $\mathrm{dlog}\, g= \lambda'-\lambda$.
\end{itemize}

The last two items define the homomorphism categories. It is convenient to
close them first under descent. This way, we obtain the prestack 
$\Grbctriv$ of trivial 
bundle gerbes with connection where the bicategories
$\Grbctriv(M)$ are defined by:

\begin{itemize}
\item
An object is a 2-form $\omega \in \Omega^2(M)$, called a
\emph{trivial bundle gerbe with connection} and denoted by 
$\mathcal{I}_{\omega}$.

\item
A 1-morphism $\mathcal{I}_\omega \to \mathcal{I}_\omega'$
is a $\Ue$ bundle $L$ with connection of curvature $\omega'-\omega$.

\item
A 2-morphism $\phi: L \to L'$
is a morphism of bundles with connection.
\end{itemize}

There is also a natural pullback operation along maps, induced by pullback 
on differential forms and pullback on $U(1)$-bundles. One checks that
$\Grbctriv$ is a prestack. By theorem \ref{stackifizierung}, the plus
construction yields a stack 
$$\Grbc:=\big(\Grbctriv\big)^+$$
on $\Man$ and even a stack on 
the category of Lie groupoids. In particular, definition \ref{defEquivariant}
provides a natural notion of an equivariant gerbe. Theorem \ref{2.13} 
now implies:

\begin{corollary} ~\\
For an equivalence $F: \Gamma \to \Lambda$ of Lie groupoids, 
the pullback functor
$$ F^*: \Grb(\Lambda) \to \Grb(\Gamma) \qquad\qquad \Grbc(\Lambda) \to \Grbc(\Gamma)$$
is an equivalence of bicategories. In particular, for a free, proper and discontinuous action of a Lie group $G$ on a smooth manifold $M$ we have
the following equivalences of bicategories
$$ \Grb_G(M) \cong \Grb(M/G) \qquad\text{respectively}\qquad  \Grbc_G(M) \cong \Grbc(M/G)
\,\, .$$
\end{corollary}

We compare the stack $\Grbc$ with objects introduced in the literature.
An object in $\Grbc(M)$ consists by definition of a covering $Y \to M$ and an object 
$G$ in $\Desc_{\Grbctriv}(Y)$. Spelling out the data explicitly, one verifies that objects
are just bundle gerbes in the sense of \cite{murray} and \cite{stev2000}. 
For the special case of an open cover $Y := \bigsqcup U_i$, an object in 
$\Desc_{\Grbctriv}(Y)$ is an Chaterjee-Hitchin gerbe, see \cite{chatterjee-gerbs}. \\

To compare different morphisms introduced in the literature,
we first need a definition:

\begin{Definition} \mbox{} \\
i) A morphism $\mathcal{A}: (Y, G) \to (Y', G')$ in $\X^+(M)$ 
is called a {\em stable isomorphism}, if it is defined on the
canonical common refinement
$$Z := Y \times_M Y'\,\, \text{.} $$
ii) A {\em stable 2-isomorphism} in $\X(M)$ between stable isomorphisms $(Z, A)$ and $(Z,A')$ 
is a morphism in $\Desc_{\X}(Z \uto M)$, i.e.\ a morphism on the canonical
common refinement $Z = Y\times_M Y'$. \\[.2em]
iii) Two objects $(Y, \calg)$ and $(Y', \calg')$ are called 
{\em stably isomorphic} if there is a stable isomorphism
$(Y, \calg) \to (Y', \calg')$.
\end{Definition}

For bundle gerbes $(Y,G)$ and $(Y',G')$, stable morphisms
are a subcategory,  
$$\Hom_{\text{Stab}}\Big((Y,G),(Y',G')\Big) \subset 
\Hom_{\Grbc}\Big((Y,G),(Y',G')\Big) \,\, . $$
We next show that the these categories are in fact
equivalent. We start with the following observation:

\begin{Lemma} \mbox{} \\
Let $Z \CTo M$ be a common refinement of $Y \CTo M$ and $Y'\CTo M$ with morphisms $s: Z \CTo Y$ and $s': Z \CTo Y'$ as in (\ref{eq:3.1}).
Then the morphism $s \times_M s': Z \to Y\times_M Y'$ induces a 
$\tau$-weak equivalence of \v{C}ech groupoids 
$$\Ce(Z) \iso \Ce(Y\times_M Y')\text{.}$$
\end{Lemma}
\begin{Proof}
Spelling out the definition of $\tau$-essential surjectivity for the Lie functor
$\Ce(Z) \to \Ce(Y\times_M Y')$, we see that we have to 
show that the smooth map
$$ Z\times_{Y\times_M Y'} \left( Y\times_M Y'\times_MY \times_M Y'\right)
\equiv Z\times_M Y\times_M Y' \to Y\times_M Y' $$
is in $\tau$. This follows at once from the pullback diagram:
$$\xymatrix{ Z\times_M Y\times_M Y' \ar[r] \ar@{->>}[d]& Z \ar@{->>}[d] \\
Y\times_M Y' \ar[r] & M 
} $$

It remains to show that the Lie functor is fully faithful. From example 
\ref{bsp:cech} we know that the vertical morphisms in the diagram
$$\xymatrix{
\Ce(Z) \ar[rr]\ar[rd] & & \Ce(Y\times_M Y') \ar[ld] \\
& M &
}$$
are $\tau$-weak equivalences, and thus in particular fully faithful. Elementary 
properties of pullback diagrams then imply that the horizontal morphism is fully 
faithful as well. 
\end{Proof}

Hence the induced
morphism of Lie groupoids
$$  \Ce(Z)  \to \Ce(Y \times_M Y')$$ 
is fully faithful and $\tau$-essentially surjective. Since $\X$ is a
prestack, we deduce from the first assertion of theorem \ref{2.13}

\begin{Proposition} \label{stabil} \mbox{}\\
For any two objects $\mathcal{O} = (Y, \calg)$ and $\mathcal{O}' = (Y', \calg')$
in $\X^+(M)$, the 1-category $\Hom(\mathcal{O},\mathcal{O}')$
is equivalent to the subcategory of stable isomorphisms and stable
2-isomorphisms. 

In particular, two objects are isomorphic in $\X^+(M)$, if and only
if they are stably isomorphic.
\end{Proposition}  

\begin{Remark} \mbox{} \\[-1.6em ]
\begin{itemize}
\item 
Stable isomorphisms have been introduced in \cite{stev2000,murray2001bgs};
proposition \ref{stabil} shows that our bicategory is equivalent to the one
in that paper. With our definition of morphisms, composition has a much
simpler structure.
\item
In \cite{waldorf2007more} a further different choice of common refinement was made.
The bicategory introduced in \cite{waldorf2007more} has as morphisms categories
that are contained in our morphism categories and contain the morphism 
categories of \cite{murray2001bgs}. Hence all three bicategories
are equivalent.
\end{itemize}
\end{Remark}

\subsection{Jandl gerbes}\label{ssec:jandl}

It is instructive to apply the same reasoning to the  
construction of Jandl gerbes. The slightly less general notion of gerbes 
with a Jandl structure has  
been introduced in \cite{ssw} to obtain a notion of surface holonomy
for unoriented surfaces. In this subsection, we introduce the more general 
notion of a Jandl gerbe. To this end, we follow the general pattern from section 
\ref{gerbes} and first define Jandl bundles:

\begin{Definition} ~\\
A {\em Jandl bundle} over $M$ is a pair, consisting of a $U(1)$-bundle $P$ 
with connection over 
$M$ and a smooth smooth map $\sigma: M \to \mathbb{Z}/2 = \{1,-1\}$. 
Morphisms of Jandl bundles $(P,\sigma) \to (Q,\mu)$ only exist
if $\sigma = \mu$. In this case they are morphisms $P\to Q$ of bundles with connection. We denote the category of Jandl bundles by 
$\JBunc(M)$
\end{Definition}

We need the covariant functor 
$$ (?)^{-1}: \,\,  \Bunc(M) \to \Bunc(M) $$
which sends a bundle $P$ to its dual bundle $P^*$. A morphism $f: P \to Q$ is sent to 
$\big(f^*\big)^{-1}: P^* \to Q^*$. This functor is well 
defined since all morphisms in $\Bunc(M)$ are isomorphisms. It squares to the 
identity and thus defines an $\mathbb{Z}/2$ action on the 
category $\Bunc(M)$.

Smooth maps $\sigma: M \to \mathbb{Z}/2$ are constant on 
connected components of $M$. For each such map $\sigma$, we
get a functor by letting $(?)^{-1}$ acting on each connected
component by the power given by the value of $\sigma$ on that
connected component. For each map $\sigma$ we thus have a
functor 
$$ (?)^{\sigma}: \,\, \Bunc(M) \to \Bunc(M)\text{.} $$ 

For our construction, we need a monoidal category of morphisms of 
trivial objects.
Hence we endow $\JBunc(M)$ with a monoidal structure;
$$(P,\sigma) \otimes (Q,\mu) := (P \otimes Q^\sigma,\sigma \mu)
\,\,\, \text{.}$$
Now we are ready to 
define the prestack $\JGrbctriv$ of trivial Jandl gerbes. 
Again the guiding principle is the definition of holonomies,
this time for unoriented surfaces (for more details, see
appendix \ref{sapp:sholunor}).
\begin{itemize}
\item
An object is a 2-form $\omega \in \Omega^2(M)$, called a
\emph{trivial Jandl gerbe with connection} and denoted by 
$\mathcal{I}_{\omega}$.

\item
A 1-morphism $\mathcal{I}_\omega \to \mathcal{I}_{\omega'}$
is a Jandl bundle $(P,\sigma)$ of curvature 
$\mathrm{curv}P=\sigma \cdot \omega'-\omega$.

\item
A 2-morphism $\phi: (P,\sigma) \to (Q,\mu)$
is a morphism of Jandl bundles with connection.
\end{itemize}
Composition of morphisms is defined as the tensor product of 
Jandl bundles. It is easy to see that $\JGrbctriv$ is a prestack. 
We define Jandl gerbes by applying the plus construction:
$$ \JGrbc := \big(\JGrbctriv\big)^+.$$
By theorem \ref{stackifizierung}, this defines a stack.

\begin{Remark} \mbox{} \\[-1.8em]
\begin{enumerate}
\item We relegate the discussion of the relation between Jandl gerbes and 
gerbes with a Jandl structure introduced in \cite{ssw} to appendix 
\ref{sapp:sholunor}, see proposition \ref{vergleichssw}. In the same
appendix, we discuss holonomy for unoriented surfaces.

\item
In terms of descent data, we can describe a Jandl gerbe on $M$ by a 
cover $Y\uto M$, a two-form $\omega\in\Omega^2(Y)$, a Jandl
bundle $(P,\sigma)$ on $Y^{[2]}$ such that
$\sigma\partial_1^*\omega-\partial_0^*\omega= \curv(P)$ and a 2-morphism 
$$\mu: \partial_2^* (P,\sigma)\otimes \partial_0^*(P,\sigma) 
\Rightarrow \partial_1^* (P,\sigma) \,\,  $$
of Jandl bundles on $Y^{[3]}$. The definition of morphisms of Jandl
bundles implies that such a morphism only exists, if the identity
\be\label{eq:cocylcecover}
\partial_2^*\sigma \cdot \partial_0^*\sigma=\partial_1^*\sigma
\ee
holds. Under this condition, the data on $Y^{[3]}$ reduce to a morphism of 
$U(1)$-bundles
$$ \mu: \,\, \partial_2^*P\otimes\partial_0^* P\Rightarrow \partial_1^* P $$
that obeys the same associativity condition on $Y^{[4]}$ as ordinary gerbes.

\item
Both trivial Jandl gerbes and trivial bundle gerbes
are given, as objects, by 2-forms; hence they are locally
the same. The crucial difference between Jandl gerbes and 
bundle gerbes is the fact that there are more 1-morphisms 
between Jandl gerbes: apart from the morphisms $(P,1)$,
we also have ``odd'' morphisms $(P,-1)$. 
\end{enumerate}
\end{Remark}

We have the inclusion $j: \Bun(M) \to \JBun(M)$ where we 
identify a bundle $P \in \Bun(M)$ with a Jandl bundle
$(P,1) \in \JBun(M)$. Here $1:M\to \mathbb{Z}/2$ is the constant function
to the neutral element. The category $\Bun(M)$ is thus a full subcategory
of $\JBun(M)$. The inclusion functor is clearly monoidal and thus yields an inclusion $\Grbctriv(M) \to \JGrbctriv(M)$ of bicategories. Finally this induces an inclusion functor
$$ \calj: \Grbc(M) \to \JGrbc(M) \,\,\, \text{.}$$
In terms of descent data, the functor $\calj$ maps
$$ (Y,\omega,P,\mu) \mapsto (Y,\omega, (P,1), \mu) \,\, . $$
The inclusion functor $\calj$ is faithful, but neither full nor essentially 
surjective. Hence we have to understand its essential image.\\

Given a Jandl bundle $(P,\sigma)$, we can forget $P$ and just
keep the smooth map $\sigma$. Since morphisms in $\JBun(M)$ preserve $\sigma$ by definition, this yields a functor
\begin{equation}\label{assignment}
 o: \JBunc(M) \to C^\infty(M,\mathbb{Z}/2) \,\,\,
\end{equation}
where the category on the right hand side has $\mathbb{Z}/2$-valued
smooth functions as objects and only identities as morphisms.

The functor $o$ is monoidal, i.e.\ $(P,\sigma) \otimes (Q,\mu) \mapsto \sigma \cdot \mu$. 
We denote the category of $\mathbb{Z}/2$ bundles on $M$ by $\Bun_{\mathbb{Z}/2}(M)$. It contains the full subcategory
$ \Buntriv_{\mathbb{Z}/2}(M)$
of trivial $\mathbb{Z}/2$-bundles:
\begin{itemize}
\item
The category $ \Buntriv_{\mathbb{Z}/2}(M)$
has exactly one object, the trivial $\mathbb{Z}/2$ bundle 
$M \times \mathbb{Z}/2\to M$.
\item
The endomorphisms of $M \times \mathbb{Z}/2$  are given by elements in $C^\infty(M,\mathbb{Z}/2)$.
\item 
Composition of endomorphisms is pointwise multiplication of smooth maps $M \to \mathbb{Z}/2$
\end{itemize}

Together with this observation the functor \eqref{assignment} yields a functor
$$ \JGrbctriv(M) \to \Buntriv_{\mathbb{Z}/2}(M).$$
Applying the plus construction, that functor induces a functor
$$ \calo: \JGrbc(M) \to \Bun_{\mathbb{Z}/2}(M).$$
In terms of descent data, the functor $\calo$ maps
$$  (Y,\omega, (P,\sigma), \mu) \mapsto (Y,\sigma) \,\, . $$
Equation (\ref{eq:cocylcecover}) implies that the cocycle condition holds
on $Y^{[3]}$
so that the pair $(Y,\sigma)$ indeed describes a $\mathbb{Z}/2$-bundle 
in terms of local data. For later use, we note that a section of the 
bundle $(Y,\sigma)$ is described in local data by a function 
$s:Y\to \mathbb{Z}/2$ such that the identity
$\sigma=\partial_0^*s\partial_1^* s$ holds on $Y^{[2]}$.

We are now ready for the next definition:

\begin{Definition} \mbox{} \label{newdef} \\[-1.8em]
\begin{enumerate} 
\item
We call $\calo(\calg)$ the \emph{orientation bundle} of the Jandl gerbe 
$\calg$.

\item A global section $s$ of the orientation bundle $\calo(\calg)$ is called
an {\em orientation} of the Jandl gerbe $\calg$.

\item
A morphism $\varphi: \calg\to \calg'$ of oriented Jandl gerbes is called 
{\em orientation preserving}, if the morphism $\calo(\varphi)$ of
$\mathbb{Z}/2$-covers preserves the global sections,
$ \calo(\varphi) \circ s = s'$.

\item
Together with all 2-morphisms of Jandl gerbes, we obtain the bicategory
$\JOGrbc(M)$ of oriented Jandl gerbes.
\end{enumerate}
\end{Definition}

\begin{Proposition}\label{sequence} \mbox{} \\[-1.8em]
\begin{enumerate}
\item \label{Teileins}
For any gerbe $\calg$, the induced Jandl gerbe $\calj(\calg)$ is canonically
oriented. For any morphism $\varphi:\calg\to\calg'$ of gerbes, the
induced morphism $\calj(\varphi): \calj(\calg)\to\calj(\calg')$ is
orientation preserving.

\item 
The functor $\calj$ induces an equivalence of bicategories
$$ \Grbc(M) \to \JOGrbc(M)  \,\, \text{.} $$
Hence the choice of an orientation reduces a Jandl gerbe to a gerbe.
\end{enumerate}
\end{Proposition}
\begin{Proof}

\begin{enumerate}
\item
Let $\calg$ be an ordinary gerbe with connection in $\Grbc(M)$.
By definition of the functors $\calj$ and $\calo$, the bundle 
$\calo(\calj(\calg))$ is given by the trivial $\mathbb{Z}/2$ cocycle on the 
covering of $\calg$. Hence it admits a canonical section $s_{\calg}$. 
This section is preserved by $\calo( \calj(\varphi))$ for any morphism 
$\varphi:\calg\to\calg'$ of gerbes. This shows part \ref{Teileins} of the 
claim. 

\item
By looking at the local data, we find that data and conditions of
a Jandl gerbe $(Y,\mathcal{I}_\omega, (P,\sigma),\mu)$ with
$\sigma: Y^{[2]}\to \mathbb{Z}/2$ the constant map to 1 are
precisely the local data of a gerbe with connection. 
Since the orientation bundle $(Y,1)$ of such a Jandl gerbe is trivial,
we choose the trivial section $1: Y\to \mathbb{Z}/2$ as the 
canonical orientation. Similarly, one sees that morphisms of such
Jandl gerbes preserving the canonical orientation are described
by exactly the same local data as morphisms of gerbes with connection.
The 2-morphisms between two such morphisms are the same anyway. 
Hence, the functor $\calj$ is an isomorphism  from the bicategory  
$\Grbc(M)$ to the full subbicategory of $\JOGrbc(M)$ with trivial 
map $\sigma$.

It remains to show that any oriented Jandl gerbe with connection is
isomorphic within $\JOGrbc(M)$ to an object in the full subbicategory
with trivial map $\sigma$.
To this end, we apply to a general Jandl gerbe 
$(Y,\mathcal{I}_\omega, (P,\sigma),\mu)$ with orientation 
$s: Y\to \mathbb{Z}/2$  the isomorphism
$m=(Y, (\mathrm{triv},s), \id)$. Here $\mathrm{triv}$ is the trivial
$U(1)$-bundle on $Y$. The target of this isomorphism is a trivially
oriented Jandl gerbe of the 
form $(Y,\mathcal{I}_{s\omega}, (P^{\partial_0^* s},1),\tilde\mu)$ and thus
in the full subbicategory of $\JOGrbc(M)$ described in the preceding paragraph.

\end{enumerate}
\end{Proof}

The last assertion crucially enters in the discussion of
unoriented surface holonomy in appendix \ref{sapp:sholunor}.

\subsection{Kapranov-Voevodsky 2-vector bundles}\label{ssec:kvvect}

As a further application of the plus construction, we 
investigate a version of 2-vector bundles, more precisely 2-vector
bundles modeled on the notion of Kapranov-Voevodsky 2-vector 
spaces \cite{kapranov19942}. The bicategory of complex
KV 2-vector spaces is (equivalent to) the following 
bicategory:
\begin{itemize}
\item
Objects are given by non negative integers $n \in \mathbb{N} = \{0,1,2,\ldots\}$. This is a shortcut for the category $\vect^n=\vect\times
\ldots \vect$, where we have the product of categories.
\item
1-morphisms $n \to m$ are given by $m\times n$ matrices $\big(V_{ij}\big)_{i,j}$ of complex vector spaces. This encodes an exact functor $\vect^n \to \vect^m$.
\item
2-morphisms $\big(V_{ij}\big)_{i,j} \Rightarrow \big(W_{ij}\big)_{i,j}$ are given by families $\big(\varphi_{ij}\big)_{i,j}$ of linear maps. This encodes a natural transformation between functors $\vect^n \to \vect^m$.
\end{itemize}

The 1-{\em iso}morphisms in this bicategory are exactly those
$n\times n$ square matrices $(V_{ij})$ for which the 
$n\times n$ matrix with non-negative integral entries
$(\dim_\C V_{ij})$ is invertible in the ring $M(n\times n,\mathbb{N})$ 
of matrices with integral entries. \\

Based on this bicategory we define for a smooth manifold $M$ 
the bicategories $\zvecttriv(M)$ of trivial 
Kapranov-Voevodsky 2-vector bundles: 
\begin{itemize}
\item
Objects are given by non negative integers $n \in \mathbb{N} = \{0,1,2,\ldots\}$. 
\item
1-morphisms $n \to m$ are given by $m\times n$ matrices $\big(E_{ij}\big)_{i,j}$ of complex vector bundles over $M$. 
\item
2-morphisms $\big(E_{ij}\big)_{i,j} \Rightarrow \big(F_{ij}\big)_{i,j}$ are given by families $\phi_{ij}: E_{ij} \to F_{ij}$ of vector bundle morphisms. 
\end{itemize}
The pullback of vector bundles turns this into a presheaf in bicategories. Since vector bundles can be glued together, the presheaf 
$\zvecttriv$ is even a prestack. Hence we can apply the
plus construction:
$$ \zvect := \Big(\zvecttriv\Big)^+.$$
By theorem \ref{stackifizierung}, we obtain a stack
of 2-vector bundles.
Thus we have properly defined bicategories of $\zvect(M)$ of 2-vector bundles over a manifold $M$ and even over Lie groupoids and thus obtained a notion
of equivariant 2-vector bundles. \\

In \cite{baas2004two} a notion of 2-vector bundles on the basis of 
Kapranov-Voevodsky 2-vector spaces has been introduced under the name of 
charted 2-vector bundles. They are defined on ordered open covers to 
accomodate more 1-isomorphisms and thus yield a richer setting for 
2-vector bundles.

\section{Proof of theorem \ref{2.13}, part 1: Factorizing 
morphisms}\label{sec:proof1}
Sections \ref{sec:proof1}--\ref{sec:proof4} are devoted
to the proof of theorem \ref{2.13}. For this proof, 
we factor any
fully faithful and $\tau$-essentially surjective Lie functor
$F:\Gamma\to\Omega$ into two morphisms of Lie groupoids 
belonging to special classes of morphisms of Lie groupoids: 
$\tau$-surjective equivalences and
strong equivalences. We first discuss these two classes of morphisms.
\subsection{Strong equivalences}\label{ssec:intequi}
We start with the definition of strong equivalences \cite{momr}. To this end,
we introduce natural transformations of Lie groupoids: 
Consider
the free groupoid on a single morphism, the interval groupoid:
$$ \I := (\I_1 \rightrightarrows \I_0)$$
It has two objects $\I_0 := \{ a, b\}$ and the four
isomorphisms $\I_1 := \{id_a, id_b, \ell, \ell^{-1}\}$ with 
$s(\ell) = a, t(\ell) = b$. Consider two functors $F,G:\c\to \d$
for two categories $\c,\d$. 
For any category $\Gamma$, we consider the
cylinder category $\Gamma \times \I$ with the canonical
inclusion functors $i_0,i_1: \Gamma \to \Gamma \times \I$.

It is an easy observation that natural isomorphisms $\eta:F\Rightarrow\ G$
are in bijection to functors $\tilde\eta:\c\times\I\to \d$ with
$\tilde\eta \circ i_0 = F$ and $\tilde\eta \circ i_1 = G$. 
(The bijection maps $\eta_c:\ F(c)\to G(c)$ to
$\tilde\eta(\id_c\times \ell)$.)

This observation allows us to reduce smoothness conditions 
on natural transformations to smoothness conditions
on functors. Hence, we consider the interval groupoid $\I$
as a discrete Lie groupoid and obtain for any Lie groupoid
$\Gamma$ the structure of a Lie groupoid on the
cylinder groupoid $\Gamma \times \I$.

\begin{Definition}\mbox{} \label{def:5.1.1} \\[-1.8em]
\begin{enumerate}
\item A {\em Lie transformation} $\eta$ between two Lie functors
$F,G: \Gamma \to \Omega$ is a Lie functor
$\eta: \Gamma \times \I \to \Omega$ with $\eta \circ i_0 = F$
and $\eta \circ i_1 = G$. 
\item Two Lie functors $F$ and $G$ are called {\em naturally 
isomorphic}, $F \simeq G$, if there exists a Lie transformation between 
$F$ and $G$.
\item A Lie functor $F : \Gamma \to \Omega$ is called an
{\em strong equivalence}, if there exists a Lie functor
$G: \Omega \to \Gamma$ such that 
$G \circ F \simeq \id_{\Gamma}$ and $F \circ G \simeq \id_{\Omega}$.
\end{enumerate}
\end{Definition}

We need the following characterization of strong
equivalences, which is completely analogous to a well-known
statement from category theory:

\begin{Proposition}\label{prop:intern} \mbox{} \\
A Lie functor $F: \Gamma \to \Omega$ is an strong equivalence
if and only if it is fully faithful and split essential surjective. The latter means that the map in definition \ref{volltreu}.2
$$ \Gamma_0 \times_{\Omega_0}\Omega_1 \to \Omega_0$$
induced by the target map has a section.
\end{Proposition}
\begin{Proof}
The proof is roughly the same as in classical category theory c.f. \cite{kassel1995quantum} Prop. XI.1.5. We only have to write down everything in diagrams, e.g. the condition fully faithful in terms of pullback diagram as in definition \ref{volltreu}. Note that the proof in \cite{kassel1995quantum} needs the axiom of choice;
in our context, we need a section of
the map $ \Gamma_0 \times_{\Omega_0}\Omega_1 \to \Omega_0$.
\end{Proof}

\begin{Lemma}\label{Lemmainterne}\mbox{} \\
If a Lie functor $F: \Gamma \to \Omega$ admits a fully
faithful retract, i.e.\ a fully faithful left inverse,
it is an strong equivalence.
\end{Lemma}
\begin{Proof}
Let $P$ be the fully faithful left inverse of $P$, hence
$$ P \circ F = \id_{\Gamma}\,\, \,\text{ .} $$
It remains to find a Lie transformation
$$ \eta: F \circ P \Longrightarrow \id_{\Omega}\,\, \text{ . } $$
Since the functor  $P$ is fully faithful, the diagram
$$
\xymatrix{
\Omega_1 \ar[rr]^{P_1}\ar[d]_{(s,t)} & &
\Gamma_1 \ar[d]^{(s,t)} \\
\Omega_0 \times \Omega_0 \ar[rr]^{P_0 \times P_0} &&
\Gamma_0 \times \Gamma_0
}
$$
is by definition \ref{volltreu} a pullback diagram. Define
$\eta: \Omega_0 \to \Omega_1 \cong \Omega_0 \times_{\Gamma_0} \Gamma_1 \times_{\Gamma_0} \Omega_0$ by
$$\eta(\omega) = \left(F_0P_0(\omega)\,,\id_{P_0(\omega)},\,\omega\right)
\,\,\, \text{ . }$$

The identities $P_0(w) = s(\id_{P_0(\omega)})$ and 
$t(\id_{P_0(\omega)}) = P_0(\omega) = P_0F_0P_0(\omega)$ imply that
this is well-defined; one also checks naturality. The two
identities

$$ s \eta(\omega) = F_0 P_0(\omega) \quad \text{ and }\quad t \eta(\omega) = \omega$$
imply that $\eta$ is indeed a Lie transformation 
$F\circ P\ \Rightarrow \id_\Omega$. One verifies that it has also
the other properties we were looking for.
\end{Proof}

\subsection{\texorpdfstring{$\tau$}{t}-surjective equivalences}\label{ssec:surjequi}

For any choice of topology $\tau$, we  introduce
the notion of $\tau$-surjective equivalence. This is called hypercover in 
\cite{zhu}.  In contrast to $\tau$-weak equivalences, $\tau$-surjective 
equivalences are required  
to be $\tau$-surjective, rather than only $\tau$-essentially 
surjective, as in definition \ref{volltreu}.

\begin{Definition}\mbox{} \\
A $\tau$-{\em surjective equivalence} is a fully faithful Lie functor
$F:\Lambda\to\Gamma$ such that $F_0:\Lambda_0\to \Gamma_0$
is a morphism in $\tau$.
\end{Definition}

\begin{Proposition}\mbox{} \label{4.4} \\
Let $F:\Lambda\to\Gamma$ be a fully faithful Lie functor
and $F_\bullet: \Lambda_\bullet\to \Gamma_\bullet$ the
associated simplicial map. Then $F$ is a $\tau$-surjective 
equivalence, if and only if all maps
$F_i:\Lambda_i\to \Gamma_i$ are in $\tau$.
\end{Proposition}

The proof is based on 

\begin{Lemma}\mbox{}\label{lemma:5.6} \\
For any two $\tau$-covers $\pi:\,\,  Y \uto M$ and 
$\pi': \,\,Y' \uto M'$ in 
$\Man$, the product 
$\pi \times \pi': Y \times Y' \to M \times M'$ 
is in $\tau$ as well.
\end{Lemma}

\begin{Proof}
Writing $\pi \times \pi' = (\pi \times id) \circ (id \times \pi')$
and using the fact that the composition of $\tau$-covers is a 
$\tau$-cover, we can assume that $\pi' = id : M'\uto M'$.
The assertion then follows from the observation that the
diagram 
$$ \xymatrix{
Y \times M' \ar[r]\ar[d] & Y\ar@{->>}^\pi[d] \\
M \times M' \ar[r] & M
}
$$
is a pullback diagram
and that $\tau$ is closed under pullbacks.
\end{Proof}

\begin{Proof} of proposition \ref{4.4}.
Since $F$\ is fully faithful, all diagrams
$$
\xymatrix{
\Lambda_n \ar[rrr]^{F_n}\ar[d] & & & \Gamma_n \ar[d] \\
{\underbrace{\Lambda_0 \times \cdots \times\Lambda_0}_{n+1}} \ar[rrr]^{F_0 \times \cdots \times F_0}
& & &
{\underbrace{\Gamma_0 \times \cdots \times\Gamma_0}_{n+1}}
}
$$
are pullback diagrams. Then $F_n$ is a $\tau$-cover since 
$F_0 \times \cdots \times F_0$ is, by lemma \ref{lemma:5.6} a $\tau$-cover.
\end{Proof}

\subsection{Factorization}\label{ssec:factorization}
\begin{Proposition}[Factorization of Lie functors]
\label{FakLemma}\mbox{} \\
Let $\Gamma$ and $\Omega$ be Lie groupoids.
Every fully faithful and $\tau$-essentially surjective Lie
functor  $F: \Gamma \to \Omega$ factors as
$$
\xymatrix{
& \Lambda \ar[rd]^H & \\
\Gamma  \ar[ru]^G \ar[rr]_F
& 
&
\Omega
}
$$
where $H$ is a $\tau$-surjective equivalence and $G$\ an strong equivalence.
\end{Proposition}

\begin{Proof}
We ensure the surjectivity of $H$ by defining
$$\Lambda_0 := \Gamma_0~_{F_0}\!\times_s \Omega_1.$$
Then $H_0 : \Lambda_0 \to \Omega_0$ is given on objects by 
the target map of $\Omega$. This is a $\tau$-covering by the
definition of $\tau$-essential surjectivity. On objects,
we define $G_0: \Gamma_0 \to \Lambda_0$ by
$\gamma \mapsto \left(\gamma,\id_{F_0(\gamma)}\right)$. This
gives the commutative diagram 
$$
\xymatrix{
& \Lambda_0 \ar[rd]^{H_0} & \\
\Gamma_0  \ar[ru]^{G_0} \ar[rr]_{F_0}
& 
&
\Omega_{0}
}
$$
on the level of objects. We combine the maps in the diagram
$$
\xymatrix{
\Gamma_1 \ar[d]_{(s,t)} \ar[rrrr]^{F_1} && & &\Omega_1 \ar[d]^{(s,t)} \\
\Gamma_0 \times \Gamma_0\ar[rr]^{G_0 \times G_0} & &
\Lambda_0 \times \Lambda_0 \ar[rr]^{H_0 \times H_0}& &
\Omega_0 \times \Omega_0
}
$$
which is a pull back diagram by definition \ref{volltreu}, since $F$ is fully faithful.
To define the Lie functor $H$ such that it is fully faithful,
we have to define $\Lambda_1$ as the pull back of the right
half of the diagram, i.e.\
$\Lambda_1 := \Lambda_0 \times_{\Omega_0} \Omega_1 \times_{\Omega_0} \Lambda_0$.
The universal property of pull backs yields a diagram

\begin{equation}\label{fac_diagramm}
\xymatrix{
\Gamma_1 \ar[d]_{(s,t)} \ar[rr]^{G_1} && \Lambda_1\ar[rr]^{H_1} \ar[d]^{(s,t)} & &\Omega_1 \ar[d]^{(s,t)} \\
\Gamma_0 \times \Gamma_0\ar[rr]^{G_0 \times G_0} & &
\Lambda_0 \times \Lambda_0 \ar[rr]^{H_0 \times H_0}& &
\Omega_0 \times \Omega_0
}
\end{equation}

in which all squares are pullbacks.
The groupoid structure on  
$\Omega = (\Omega_1 \gto\Omega_0)$ 
induces a groupoid structure on 
$\Lambda = (\Lambda_1 \gto \Lambda_0)$ in such a way that
$G$ and $H$ become Lie functors.

By construction of this factorization, $H$ is a $\tau$-surjective equivalence.\ It remains to be shown that $G$ is an
strong equivalence. According to proposition \ref{prop:intern},
it suffices to show that $G$ is fully faithful and split essential surjective. The left diagram in $\eqref{fac_diagramm}$ is a pullback diagram. Hence $G$ is fully faithful. It remains to give a section of the map
\begin{equation}\label{Schnitt}
\Gamma_0 \times_{\Lambda_0} \Lambda_1 \to \Lambda_0
\end{equation}
Since we have defined $\Lambda_1 = \Lambda_0 \times_{\Omega_0} \Omega_1 \times_{\Omega_0} \Lambda_0$, we have 
$$ \Gamma_0 \times_{\Lambda_0} \Lambda_1 \cong \Gamma_0 \times_{\Omega_0} \Omega_1 \times_{\Omega_0} \Lambda_0.$$
Thus a section of \eqref{Schnitt} is given by three maps
$$ \Lambda_0 \to \Gamma_0 \qquad \Lambda_0 \to \Omega_1 \qquad \Lambda_0 \to \Lambda_0$$
that agree on $\Omega_0$, when composed with the source and the 
target map of $\Omega_0$, respectively. By definition 
$\Lambda_0 = \Gamma_0~_{F_0}\!\times_s \Omega_1$, and 
we can define the three maps by projection to the first factor, projection to the second factor and the identity.
\end{Proof}

The factorization lemma allows to isolate the violation
of $\tau$-surjectivity in an strong equivalence and to work with
$\tau$-surjective equivalences rather than only $\tau$-essentially
surjective equivalences.
Hence it suffices to prove theorem \ref{2.13} for
$\tau$-surjective equivalences and for strong equivalences.
This will be done in sections \ref{sec:proof2}
and \ref{sec:proof4}, respectively.

\section{Proof of theorem \ref{2.13}, part 2: Sheaves and strong
equivalences}\label{sec:proof2}

\begin{Lemma} \mbox{} \label{lem:6.1}\\
Let $\X$ be a presheaf that preserves products, cf.\ equation
 \eqref{disjoint_union}.
Let $\Gamma$ be a Lie groupoid and $D$ be a discrete Lie 
groupoids i.e.\ $D_0$ and $D_1$ are discrete manifolds. Then $D$ can also be  regarded as a bicategory and we have natural equivalences
$$ \X(\Gamma \times D) \cong \Big[D,\X(\Gamma)\Big]$$
where $[D,\X(\Gamma)]$ denotes the bicategory of functors $D \to \X(\Gamma)$.
\end{Lemma}

\begin{Proof}
The claim is merely a consequence the requirement 
\eqref{disjoint_union} that $\X$ preserves products: In the case that $\Gamma$ is  a manifold $M$ considered as a Lie groupoid and $D$ a set $I$ considered as a discrete groupoid we have $M \times I = \bigsqcup_{i \in I}M$. Thus the left hand is equal to $\X(\bigsqcup M)$ and the right hand side to $\X(M)^I = \prod_{i \in I} \X(M)$. In this case \eqref{disjoint_union} directly implies the equivalence. 

In the case of a general Lie groupoid, the product $\Gamma \times D$ decomposes levelwise into a disjoint union. Using this fact and explicitly spelling out $\X(\Gamma \times D)$ and $[D,\X(\Gamma)]$ according to 
definition, \ref{defEquivariant} it is straightforward to 
see that the two bicategories are equivalent.
\end{Proof}

\begin{Proposition} \mbox{} \\
Let $\X$ be a presheaf in bicategories.
Any Lie transformation $\eta:F\Rightarrow G$ of Lie functors
$F,G:\ \Gamma\to\Omega$ induces a natural isomorphism
of the pullback functors
$F^*, G^* : \X(\Omega) \to \X(\Gamma)$.
\end{Proposition}
\begin{Proof}
Recall from definition \ref{def:5.1.1} that the Lie transformation $\eta$ is by definition a Lie functor 
$$\Gamma  \times I \to \Omega \,\,\, , $$
where $I$ is the interval groupoid. Applying the presheaf $\X$ 
to this functor yields a functor
$$\X(\Omega) \to \X(\Gamma \times I).$$
Since $I$ is discrete the preceding lemma \ref{lem:6.1} shows that this is a functor
$$\X(\Omega) \to \Big[I,\X(\Gamma)\Big].$$
That is the same as a functor 
$$\X(\Omega) \times I \to \X(\Gamma)$$
i.e. a natural isomorphism of bifunctors.
\end{Proof}

\begin{Kor}\label{intern}\mbox{} \\
For any presheaf $\X$ in bicategories, the pull back along
an strong equivalence $\Gamma \to \Omega$ induces an
equivalence $\X(\Omega) \to \X(\Gamma)$ of bicategories.
\end{Kor}

\section{Proof of theorem \ref{2.13}, part 3: Equivariant descent}\label{sec:proof3}

To deal with $\tau$-surjective equivalences, we need to consider simplicial
objects in the category of simplicial objects, i.e.\
bisimplicial objects. In the course of our investigations,
we obtain results about
bisimplicial objects that are of independent interest,
in particular theorem \ref{HauptsatzEquivDesc} and corollary
\ref{EqDesc} on equivariant descent.

We first generalize the definition of equivariant objects as follows:
If we evaluate a presheaf in bicategories $\X$ on a
simplicial object $\Gamma_\bullet$, we obtain the following
diagram in $\BiCat$:
$$
\xymatrix{
{\X(\Gamma_0)}
\ar@<0.3ex>[r]^-{\partial_0^*}
\ar@<-0.7ex>[r]_-{\partial_1^*} 
&
{\X(\Gamma_1)}
\ar@<0.9ex>[r]^{\partial_0^*}
\ar@<-0.1ex>[r]
\ar@<-1.1ex>[r]_{\partial_2^*}
&
{\X(\Gamma_2)}
\ar@<1.3ex>[r]^-{\partial_0^*}
\ar@<0.3ex>[r]
\ar@<-0.7ex>[r]
\ar@<-1.7ex>[r]_-{\partial_3^*}
&
{\cdots}  
}
$$
in which the cosimplicial identities are obeyed up to natural isomorphism,
\begin{eqnarray*} 
  \partial_{j}^* \partial_{i}^* &\cong& \partial_{i}^* \partial_{j-1}^* 
  \qquad\text{for } i<j \,\, .
\end{eqnarray*}
The coherence cells turn this into a weak functor $\Delta \to \BiCat$ 
from the simplicial category $\Delta$ to $\BiCat$. Such a functor 
will be called a (weak) cosimplicial bicategory.

The equivariant objects can be constructed in this framework
by selecting objects in $\X(\Gamma_0)$, 1-morphisms in
$\X(\Gamma_1)$ and so on. This leads us to the following definition:

\begin{Definition}\mbox{} \\
Given a cosimplicial bicategory
$C_\bullet$, we introduce the category 
$$
\holim_{i \in \Delta} C_i \equiv \holim\left( 
\xymatrix{
C_0
\ar@<0.3ex>[r]^-{\partial_0^*}
\ar@<-0.7ex>[r]_-{\partial_1^*} 
&
C_1
\ar@<0.9ex>[r]^{\partial_0^*}
\ar@<-0.1ex>[r]
\ar@<-1.1ex>[r]_{\partial_2^*}
&
C_2
\ar@<1.3ex>[r]^-{\partial_0^*}
\ar@<0.3ex>[r]
\ar@<-0.7ex>[r]
\ar@<-1.7ex>[r]_-{\partial_3^*}
&
{\cdots}  
}
\right)
$$
with objects given by the following data:

\begin{itemize}
\item[(O1)] An object $\mathcal{G}$ in the bicategory $C_0$;
\item[(O2)] 
A 1-isomorphism in the bicategory $C_1$;
  \begin{equation*}
  P\colon~ \partial_0^{*}\mathcal{G} \To \partial_1^{*}\mathcal{G}
  \end{equation*} 

\item[(O3)] 
A 2-isomorphism in the bicategory $C_2$;
  \begin{equation*}
  \mu \colon~ \partial_{2}^*P \oti \partial_{0}^*P \,{\Rightarrow}\, \partial_{1}^*P
  \end{equation*} 

\item[(O4)] 
A coherence condition of 2-morphisms in the bicategory
$C_3$:

  \begin{equation*}
  \partial_2^*\mu \circ (\id \oti \partial_0^* \mu) = \partial_1^*\mu \circ (\partial_3^*\mu \oti \id)\,
  \end{equation*}
\end{itemize}
Morphisms and 2-morphisms are defined as in definition \ref{2.6}.
\end{Definition}

In this notation, the extension of a prestack $\X$
to an equivariant object $\Gamma_\bullet$ described in
definition \ref{defEquivariant} is given by
\be
\X (\Gamma_\bullet) = \holim_{i \in \Delta} \X(\Gamma_i)\,\,.
\label{Nummer}
\ee
In the special case of a $\tau$-covering $Y\uto M$, we can
write the descent object as
$$\Desc_\X (Y \uto M) = \holim_{i \in \Delta} \X(Y^{[i+1]}).$$
For the constant simplicial bicategory $C_\bullet$, with 
$C_i = C$ for all $i$, one checks that
$\holim_{i \in \Delta} C_i = C$.

We next need to extend the notion of a $\tau$-covering to
a simplicial object:

\begin{Definition}\label{CoveringGruppoids}\mbox{} \\[-1.8em]
\begin{enumerate}
\item
Let $\Lambda_\bullet$  and $\Gamma_\bullet$ be simplicial manifolds and 
$\Pi_\bullet: \Lambda_\bullet \to \Gamma_\bullet$ a simplicial
map. Then $\Pi_\bullet$ is called a $\tau$-cover, if all maps
$\Pi_i : \Lambda_i \to \Gamma_i$ are $\tau$-covers.
\item
A Lie functor $\Pi: (\Lambda_1\rightrightarrows \Lambda_0) \to (\Gamma_{1} \rightrightarrows \Gamma_{0})$ 
is called a $\tau$-cover, if the associated simplicial map
$ \Pi_\bullet : \Lambda_{\bullet} \to \Gamma_{\bullet}$ 
of the nerves is a $\tau$-cover of simplicial manifolds.
\end{enumerate}
\end{Definition}

\begin{Remark}\label{bemerkung}\mbox{} \\[-1.8em]
\begin{enumerate}
\item
Proposition \ref{4.4} shows that for a $\tau$-surjective equivalence the associated simplicial map is $\tau$-cover.
\item
For any $\tau$-covering $\pi: Y \uto M$, the simplicial
map induced by the Lie functor $\Ce(Y) \to M$ is an
example of a $\tau$-cover of simplicial manifolds.
\end{enumerate}
\end{Remark}

Given a $\tau$-cover $\Pi_\bullet: \Lambda_\bullet \to \Gamma_\bullet$
of simplicial manifolds, we can construct the simplicial
manifold
$$
\Lambda^{[2]}_\bullet:=
\Lambda_\bullet \times_{\Gamma_\bullet} \Lambda_\bullet := 
\left( 
\xymatrix{
{\cdots}  
\ar@<1.3ex>[r]^-{\partial_0}
\ar@<0.3ex>[r]
\ar@<-0.7ex>[r]
\ar@<-1.7ex>[r]_-{\partial_3}
&
{\Lambda_2 \times_{\Gamma_2} \Lambda_2}
\ar@<0.9ex>[r]^{\partial_0}
\ar@<-0.1ex>[r]
\ar@<-1.1ex>[r]_{\partial_2}
&
{\Lambda_1 \times_{\Gamma_1} \Lambda_1}
\ar@<0.3ex>[r]^-{\partial_0}
\ar@<-0.7ex>[r]_-{\partial_1} 
&
{\Lambda_0 \times_{\Gamma_0} \Lambda_0}
}
\right)
$$
with obvious maps  $\partial_i$. One verifies that the
two projections 
 $\delta_0, \delta_1: {\Lambda_\bullet^{[2]}} \to \Lambda_\bullet$
are simplicial maps. Similarly, we form simplicial manifolds
$$\Lambda^{[n]}_\bullet := \underbrace{\Lambda_\bullet \times_{\Gamma_\bullet} \ldots \times_{\Gamma_\bullet} \Lambda_\bullet}_n$$ 
and simplicial maps $\delta_i: \Lambda^{[n]}_\bullet \to \Lambda^{[n-1]}_\bullet$ with $i=0,\ldots,n-1$. We thus obtain an (augmented) simplicial
object
$$
\big(\Lambda_{\bullet}\big)^{[\bullet]} ~ := ~
\left( 
\xymatrix{
{\cdots}  
\ar@<1.3ex>[r]^-{\delta_0}
\ar@<0.3ex>[r]
\ar@<-0.7ex>[r]
\ar@<-1.7ex>[r]_-{\delta_3}
&
{\Lambda^{[3]}_\bullet}
\ar@<0.9ex>[r]^{\delta_0}
\ar@<-0.1ex>[r]
\ar@<-1.1ex>[r]_{\delta_2}
&
{\Lambda^{[2]}_\bullet}
\ar@<0.3ex>[r]^-{\delta_0}
\ar@<-0.7ex>[r]_-{\delta_1} 
&
{\Lambda_\bullet}
}
\right)
\longrightarrow \Gamma_\bullet
$$
in the category of simplicial manifolds. A simplicial
object in the category of simplicial manifolds will also
be called a {\em bisimplicial manifold}. In full detail,
a bisimplicial manifold consists of the following data:
$$\xymatrix{
 & 
\cdots 
\ar@<1.3ex>[d]
\ar@<0.3ex>[d]
\ar@<-0.7ex>[d]
\ar@<-1.7ex>[d]
&
\cdots 
\ar@<1.3ex>[d]
\ar@<0.3ex>[d]
\ar@<-0.7ex>[d]
\ar@<-1.7ex>[d]
&
\cdots 
\ar@<1.3ex>[d]
\ar@<0.3ex>[d]
\ar@<-0.7ex>[d]
\ar@<-1.7ex>[d]
&& 
\cdots
\ar@<1.3ex>[d]
\ar@<0.3ex>[d]
\ar@<-0.7ex>[d]
\ar@<-1.7ex>[d]
\\
\cdots
\ar@<1.3ex>[r]
\ar@<0.3ex>[r]
\ar@<-0.7ex>[r]
\ar@<-1.7ex>[r]
& 
\Lambda^{[3]}_2
\ar@<0.9ex>[r]
\ar@<-0.1ex>[r]
\ar@<-1.1ex>[r]
\ar@<0.9ex>[d]
\ar@<-0.1ex>[d]
\ar@<-1.1ex>[d]
&
\Lambda^{[2]}_2
\ar@<0.3ex>[r]
\ar@<-0.7ex>[r]
\ar@<0.9ex>[d]
\ar@<-0.1ex>[d]
\ar@<-1.1ex>[d]
&
\Lambda_2
\ar@<0.9ex>[d]
\ar@<-0.1ex>[d]
\ar@<-1.1ex>[d]
\ar[rr]
&& 
\Gamma_2
\ar@<0.9ex>[d]
\ar@<-0.1ex>[d]
\ar@<-1.1ex>[d]
\\
\cdots
\ar@<1.3ex>[r]
\ar@<0.3ex>[r]
\ar@<-0.7ex>[r]
\ar@<-1.7ex>[r]
&
\Lambda^{[3]}_1
\ar@<0.9ex>[r]
\ar@<-0.1ex>[r]
\ar@<-1.1ex>[r]
\ar@<0.3ex>[d]
\ar@<-0.7ex>[d]
&
\Lambda^{[2]}_1
\ar@<0.3ex>[r]
\ar@<-0.7ex>[r]
\ar@<0.3ex>[d]
\ar@<-0.7ex>[d]
&
\Lambda_1
\ar@<0.3ex>[d]
\ar@<-0.7ex>[d]
\ar[rr]
&& 
\Gamma_1
\ar@<0.3ex>[d]
\ar@<-0.7ex>[d]
\\
\cdots
\ar@<1.3ex>[r]
\ar@<0.3ex>[r]
\ar@<-0.7ex>[r]
\ar@<-1.7ex>[r]
&
\Lambda^{[3]}_0
\ar@<0.9ex>[r]
\ar@<-0.1ex>[r]
\ar@<-1.1ex>[r]
&
\Lambda^{[2]}_0
\ar@<0.3ex>[r]
\ar@<-0.7ex>[r]
&
\Lambda_0
\ar[rr]
&& \Gamma_0
}
$$
The rows are, by construction, nerves of \v Cech groupoids.
This fact will enter crucially in the proof of our
main result on equivariant descent. Before turning to this,
we need the following 

\begin{Proposition}\mbox{}\label{bisimplicial} \\
Let $\X$ be a presheaf in bicategories and $\Omega_{\bullet \bullet}$
a bisimplicial manifold. Then
$$
\holim_{i \in \Delta} \holim_{j \in \Delta} \X\left(\Omega_{ij}\right) =  \holim_{j \in \Delta} \holim_{i \in \Delta} \X\left(\Omega_{ij}\right)
$$
\end{Proposition}

\begin{Proof}
We first discuss what data of the bisimplicial manifold
$$\xymatrix{
 & 
\cdots 
\ar@<1.3ex>[d]
\ar@<0.3ex>[d]
\ar@<-0.7ex>[d]
\ar@<-1.7ex>[d]
&
\cdots 
\ar@<1.3ex>[d]
\ar@<0.3ex>[d]
\ar@<-0.7ex>[d]
\ar@<-1.7ex>[d]
&
\cdots 
\ar@<1.3ex>[d]
\ar@<0.3ex>[d]
\ar@<-0.7ex>[d]
\ar@<-1.7ex>[d]
\\
\cdots
\ar@<1.3ex>[r]
\ar@<0.3ex>[r]
\ar@<-0.7ex>[r]
\ar@<-1.7ex>[r]
& 
\Omega_{22}
\ar@<0.9ex>[r]
\ar@<-0.1ex>[r]
\ar@<-1.1ex>[r]
\ar@<0.9ex>[d]
\ar@<-0.1ex>[d]
\ar@<-1.1ex>[d]
&
\Omega_{21}
\ar@<0.3ex>[r]
\ar@<-0.7ex>[r]
\ar@<0.9ex>[d]
\ar@<-0.1ex>[d]
\ar@<-1.1ex>[d]
&
\Omega_{20}
\ar@<0.9ex>[d]
\ar@<-0.1ex>[d]
\ar@<-1.1ex>[d]
\\
\cdots
\ar@<1.3ex>[r]
\ar@<0.3ex>[r]
\ar@<-0.7ex>[r]
\ar@<-1.7ex>[r]
&
\Omega_{12}
\ar@<0.9ex>[r]
\ar@<-0.1ex>[r]
\ar@<-1.1ex>[r]
\ar@<0.3ex>[d]
\ar@<-0.7ex>[d]
&
\Omega_{11}
\ar@<0.3ex>[r]
\ar@<-0.7ex>[r]
\ar@<0.3ex>[d]
\ar@<-0.7ex>[d]
&
\Omega_{10}
\ar@<0.3ex>[d]
\ar@<-0.7ex>[d]
\\
\cdots
\ar@<1.3ex>[r]
\ar@<0.3ex>[r]
\ar@<-0.7ex>[r]
\ar@<-1.7ex>[r]
&
\Omega_{02}
\ar@<0.9ex>[r]
\ar@<-0.1ex>[r]
\ar@<-1.1ex>[r]
&
\Omega_{01}
\ar@<0.3ex>[r]
\ar@<-0.7ex>[r]
&
\Omega_{00}
}
$$
enter in an object of 
$\holim_{i \in \Delta} \holim_{j \in \Delta} \X\left(\Omega_{ij}\right)$.
To this end, we denote horizontal
boundary maps by $\delta$ and vertical boundary maps by
$\partial$. Then such an object is given by
\begin{itemize}
\item An object in
$\holim_{j \in \Delta} \X\left(\Omega_{0j}\right)$ which in
turn consists of 
    \begin{itemize}
    \item An object $\calg$ on $\Omega_{00}$
    \item An isomorphism $A_{01}: \delta_0^* \calg \to \delta_1^* \calg$  
    on $\Omega_{01}$
    \item A 2-isomorphism 
    $\mu_{02}: \delta_2^* A_{01} \otimes \delta_0^* A_{01} \Rightarrow \delta_1^*A_{01}$ 
    on $\Omega_{02}$
    \item A coherence condition on  $\Omega_{03}$
    \end{itemize}    
\item A morphism $\partial_0^*(\calg, A_{0,1},\mu_{0,2}) \to \partial_1^*(\calg, A_{0,1},\mu_{02})$ in $\holim_{j \in \Delta} \X\left(\Omega_{1j}\right)$ which in turn consists of 
     \begin{itemize}
     \item An isomorphism $A_{10}: \partial_0^* \calg 
     \to \partial_1^* \calg$ on $\Omega_{10}$
     \item A 2-isomorphism
           $\mu_{11}: \partial_1^* A_{01} \otimes \delta_0^* A_{10} \Rightarrow \delta_1^* A_{10} \otimes \partial_0^* A_{01}$ on $\Omega_{11}$.
     \item A coherence condition on $\Omega_{12}$.
     \end{itemize}
\item A 2-isomorphism $\partial_2^*(A_{10},\mu_{11}) \otimes \partial_0^*(A_{10},\mu_{11}) \Rightarrow \partial_1^*(A_{10},\mu_{11})$ in $\holim_{j \in \Delta} \X\left(\Omega_{2j}\right)$:
     \begin{itemize}
     \item A 2-isomorphism $\mu_{20}: \partial_2^*A_{10} \otimes \partial_0^* A_{10} \Rightarrow \partial_1^* A_{10}$ on $\Omega_{20}$.
     \item A coherence condition on  $\Omega_{21}$.
     \end{itemize}
\item A condition on the 2-morphisms 
in $\holim_{j \in \Delta} \X\left(\Omega_{3j}\right)$ which is
just 
\begin{itemize}
     \item A coherence condition on  $\Omega_{30}$.
\end{itemize}
\end{itemize}
To summarize, we get an object 
$\calg \in \X(\Omega_{00})$ in the lower right corner of
the diagram, two isomorphisms 
$A_{01} \in \X(\Omega_{01})$, $A_{10} \in \X(\Omega_{01})$ 
on the diagonal, three 2-isomorphisms 
$\mu_{02} \in \X(\Omega_{02}), \mu_{11} \in \X(\Omega_{11}), \mu_{20} \in \X(\Omega_{20})$ on the first translate of the diagonal
and four conditions on the second translate of the
diagonal. 

For an object in $\holim_{j \in \Delta} \holim_{i \in \Delta} \X\left(\Omega_{ij}\right)$,
we get the same data, as can be seen by exchanging the roles of $i$ and
$j$. Since we interchange the roles of $\partial$ and $\delta$,
we have to replace the 2-isomorphism 
$\mu_{11}: \partial_1^* A_{01} \otimes \delta_0^* A_{10} \Rightarrow \delta_1^* A_{10} \otimes \partial_0^* A_{01}$ by its inverse.
For all other isomorphisms and conditions, the objects remain unchanged.

By analogous considerations, one also checks that the
morphisms and 2-morphisms in both bicategories coincide.
\end{Proof}

\begin{Theorem}[Equivariant descent]
\label{HauptsatzEquivDesc} \mbox{} \\
Let $\Pi: \Lambda_\bullet \to \Gamma_\bullet$ be a 
$\tau$-covering of
simplicial manifolds.
\begin{enumerate}
\item Let $\X$ be a $\tau$-stack
   on $\Man$.
Then we have the following equivalence of bicategories:
\begin{equation*}
\X\left(\Gamma_\bullet\right) \iso
\holim \left(
\xymatrix{
{\X\left(\Lambda_\bullet\right)}
\ar@<0.3ex>[r]^-{\delta_0^*}
\ar@<-0.7ex>[r]_-{\delta_1^*} 
&
{\X\left(\Lambda^{[2]}_\bullet\right)}
\ar@<0.9ex>[r]^{\delta_0^*}
\ar@<-0.1ex>[r]
\ar@<-1.1ex>[r]_{\delta_2^*}
&
{\X\left(\Lambda^{[3]}_\bullet\right)}
\ar@<1.3ex>[r]^-{\delta_0^*}
\ar@<0.3ex>[r]
\ar@<-0.7ex>[r]
\ar@<-1.7ex>[r]_-{\delta_3^*}
&
{\cdots}  
}
\right)
\end{equation*}
In other words, we have extended $\X$ to a $\tau$-stack on the
category of simplicial manifolds.
\item
If $\X$ is a $\tau$-prestack on $\Man$, this functor is still
fully faithful, i.e.\ an equivalence of the $\Hom$-categories.
\end{enumerate}
\end{Theorem}

\begin{Proof}
By definition, we have  $\X(\Gamma_\bullet) = \holim_{i \in \Delta}\X(\Gamma_i)$.
Since $\X$ is supposed to be a $\tau$-stack and since
all $\Pi_i: \Lambda_i \uto \Gamma_i$ are $\tau$-covers, we have
the following equivalence of bicategories:
$$ \X(\Gamma_i) \iso \Desc_\X(\Lambda_i \uto \Gamma_i) = \holim_{j \in \Delta} \X(\Lambda^{[j]}_i) \,\, .$$
Altogether, we have the equivalence of
bicategories
$$ \X(\Gamma_\bullet) \iso \holim_{i \in \Delta} \holim_{j \in \Delta} \X\left(\Lambda^{[j]}_i\right)$$
By proposition \ref{bisimplicial}, we can exchange the
homotopy limits and get

\begin{eqnarray*}
\holim_{i \in \Delta} \holim_{j \in \Delta} \X\left(\Lambda^{[j]}_i\right) &=& 
\holim_{j \in \Delta} \holim_{i \in \Delta} \X\left(\Lambda^{[j]}_i\right) \\
&\stackrel{(\ref{Nummer})}=& \holim_{j \in \Delta} \X\left(\Lambda^{[j]}_\bullet\right) 
\end{eqnarray*}
and thus the assertion for stacks.\ The assertion in the
case when $\X$ is a prestack follows by an analogous
argument.

\end{Proof}

By restriction, we obtain a $\tau$-stack on 
the full subcategory of Lie groupoids. By a further 
restriction, we get a $\tau$-stack
on the full subcategory of $G$-manifolds.
For convenience, we state our result in the special case
of $G$-manifolds:

\begin{Kor}\mbox{}\label{EqDesc} \\
Let $M$ be a $G$-manifold and $\{U_i\}_{i\in I}$ be a $G$-invariant
covering. Denote, as usual $Y:=\sqcup_{i\in I} U_i$.
Then we have:
\begin{equation*}
\X_G\left(M\right) \iso
\holim \left(
\xymatrix{
{\X_G\left(Y\right)}
\ar@<0.3ex>[r]^-{\delta_0^*}
\ar@<-0.7ex>[r]_-{\delta_1^*} 
&
{\X_G\left(Y^{[2]}\right)}
\ar@<0.9ex>[r]^{\delta_0^*}
\ar@<-0.1ex>[r]
\ar@<-1.1ex>[r]_{\delta_2^*}
&
{\X_G\left(Y^{[3]}\right)}
\ar@<1.3ex>[r]^-{\delta_0^*}
\ar@<0.3ex>[r]
\ar@<-0.7ex>[r]
\ar@<-1.7ex>[r]_-{\delta_3^*}
&
{\cdots}  
}
\right)
\end{equation*}
\end{Kor}

\section{Proof of theorem \ref{2.13}, part 4: Sheaves and \texorpdfstring{$\tau$}{t}-surjective equivalences}\label{sec:proof4}

We are now ready to prove theorem \ref{2.13} in the special case
of $\tau$-surjective equivalences. This actually finishes 
the proof of theorem \ref{2.13}, since by
the factorization lemma \ref{FakLemma} we have to consider only
the two cases of $\tau$-surjective equivalences and strong
equivalences.\ The latter case has already been settled with
corollary \ref{intern}.
We start with the following 

\begin{Lemma} \label{lemma2}\mbox{} \\
Let $F: \Gamma \to \Lambda$ be a $\tau$-surjective equivalence of Lie
groupoids. By remark \ref{bemerkung}.1, the functor $F$ induces 
a $\tau$-cover of Lie groupoids.
\begin{enumerate}
\item[(i)]
For any $n$, we have a canonical functor
$M^n: \Gamma^{[n]} \to \Lambda $ which is given by
arbitrary compositions in the augmented simplicial
manifold
$$
\xymatrix{
{\cdots}  
\ar@<1.3ex>[r]^-{\delta_0}
\ar@<0.3ex>[r]
\ar@<-0.7ex>[r]
\ar@<-1.7ex>[r]_-{\delta_3}
&
{\Gamma^{[3]}}
\ar@<0.9ex>[r]^{\delta_0}
\ar@<-0.1ex>[r]
\ar@<-1.1ex>[r]_{\delta_2}
&
{\Gamma^{[2]}}
\ar@<0.3ex>[r]^-{\delta_0}
\ar@<-0.7ex>[r]_-{\delta_1} 
&
{\Gamma} \ar[r]^F
&
\Lambda
}
$$
Then the functor $M^n$ is a $\tau$-surjective equivalence.

\item[(ii)]
The diagonal functors $\Gamma \to \Gamma^{[n]}$ are strong
equivalences.
\end{enumerate}
\end{Lemma}

\begin{Proof}
(i) As compositions of $\tau$-coverings, all functors
$M^n$ are $\tau$-coverings. The functor 
$F: \Gamma \to \Lambda$ is in particular fully faithful.
Hence, 
$$\Gamma_1 = \Gamma_0 \times_{\Lambda_0} \Lambda_1 
\times_{\Lambda_0} \Gamma_0\,\, \text{.}$$
We now calculate 
\begin{eqnarray*}
\Gamma_1^{[n]} &=& \Gamma_1 \times_{\Lambda_1} \cdots \times_{\Lambda_1} \Gamma_1 \\ 
&=&
\Big(\Gamma_0 \times_{\Lambda_0} \Lambda_1 \times_{\Lambda_0} \Gamma_0\Big)
\times_{\Lambda_1} \cdots \times_{\Lambda_1}
\Big(\Gamma_0 \times_{\Lambda_0} \Lambda_1 \times_{\Lambda_0} \Gamma_0\Big) 
\end{eqnarray*}
and find by reordering that this equals
$$ 
\Big(\Gamma_0 \times_{\Lambda_1}\cdots \times_{\Lambda_1}\Gamma_0\Big) 
\times_{\Lambda_0} \Lambda_1 \times_{\Lambda_0}
\Big(\Gamma_0 \times_{\Lambda_1}\cdots \times_{\Lambda_1}\Gamma_0\Big) 
\,\, . $$
Hence the diagram
\begin{equation}\label{pullbackdiag}
\xymatrix{
\Gamma_1^{[n]} \ar[rr]^{M^n_1}\ar[d] && \Lambda_1\ar[d] \\
\Gamma_0^{[n]} \times \Gamma_0^{[n]} \ar[rr]^{M^n_0} && \Lambda_0 \times \Lambda_0
}
\end{equation}
is a pullback diagram and thus the functor $M^n$ is fully faithful.

\noindent
(ii) Take any of the $n$ possible projection functors
$P^n: \Gamma^{[n]} \to \Gamma$ and consider the diagram
$$
\xymatrix{
\Gamma_1^{[n]} \ar[rr]^{P^n_1}\ar[d] && \Gamma_1 \ar[rr]^{F_1}\ar[d] && \Lambda_1\ar[d] \\
\Gamma_0^{[n]} \times \Gamma_0^{[n]} \ar[rr]^{P^n_0 \times P^n_0} && \Gamma_0 \times \Gamma_0 \ar[rr]^{F_0 \times F_0} && \Lambda_0 \times \Lambda_0
}
$$
The right diagram is by our assumptions on $F$ a pullback diagram. The external
diagram is just the diagram (\ref{pullbackdiag}) considered
in part (i)\ of the lemma and thus a pullback diagramm, as well.
Hence also the left part of the diagram is a pullback diagram
and thus the functor $P^n$ is fully faithful.
The functor $P^n$ is a left inverse of of the diagonal functor
$\Lambda \to \Lambda^{[n]}$. Lemma \ref{Lemmainterne} now
implies that the diagonal functors are strong equivalences.

\end{Proof}
 
\begin{Proposition}\label{hauptsatz3}\mbox{} \\
Let $\X$ be a presheaf in bicategories and $\Gamma \to \Lambda$
be a $\tau$-surjective equivalence. Then we have the following
equivalences of bicategories
\begin{eqnarray*}
\X(\Gamma_\bullet) &\cong& 
\holim \left(
\xymatrix{
{\X\left(\Gamma_\bullet\right)}
\ar@<0.3ex>[r]^-{\delta_0^*}
\ar@<-0.7ex>[r]_-{\delta_1^*} 
&
{\X\left(\Gamma^{[2]}_\bullet\right)}
\ar@<0.9ex>[r]^{\delta_0^*}
\ar@<-0.1ex>[r]
\ar@<-1.1ex>[r]_{\delta_2^*}
&
{\X\left(\Gamma^{[3]}_\bullet\right)}
\ar@<1.3ex>[r]^-{\delta_0^*}
\ar@<0.3ex>[r]
\ar@<-0.7ex>[r]
\ar@<-1.7ex>[r]_-{\delta_3^*}
&
{\cdots}  
}
\right) \\
&\cong&
\holim \left(
\xymatrix{
{\Desc_\X(\Gamma_0 \uto \Lambda_0)}
\ar@<0.3ex>[r]^-{\partial_0^*}
\ar@<-0.7ex>[r]_-{\partial_1^*} 
&
{\Desc_\X(\Gamma_1 \uto \Lambda_1)}
\ar@<0.9ex>[r]^{\partial_0^*}
\ar@<-0.1ex>[r]
\ar@<-1.1ex>[r]_{\partial_2^*}
&
{\Desc_\X(\Gamma_2 \uto \Lambda_2)}
\ar@<1.3ex>[r]^-{\partial_0^*}
\ar@<0.3ex>[r]
\ar@<-0.7ex>[r]
\ar@<-1.7ex>[r]_-{\partial_3^*}
&
{\cdots}  
}
\right)
\end{eqnarray*}
\end{Proposition}
 
\begin{Proof}
The diagonal functors $\Gamma \to \Gamma^{n}$ give a morphism
of simplicial manifolds 
$$
\xymatrix{
{\cdots}  
\ar@<1.3ex>[r]^-{\delta_0}
\ar@<0.3ex>[r]
\ar@<-0.7ex>[r]
\ar@<-1.7ex>[r]_-{\delta_3}
&
{\Gamma^{[2]}}
\ar@<0.9ex>[r]^{\delta_0}
\ar@<-0.1ex>[r]
\ar@<-1.1ex>[r]_{\delta_2}
&
{\Gamma^{[1]}}
\ar@<0.3ex>[r]^-{\delta_0}
\ar@<-0.7ex>[r]_-{\delta_1} 
&
{\Gamma} \ar[r]^F
&
\Lambda \ar@{=}[d]
\\
{\cdots}  
\ar@<1.3ex>[r]^-{\delta_0}
\ar@<0.3ex>[r]
\ar@<-0.7ex>[r]
\ar@<-1.7ex>[r]_-{\delta_3}
&
{\Gamma}
\ar@<0.9ex>[r]^{\delta_0}
\ar@<-0.1ex>[r]
\ar@<-1.1ex>[r]_{\delta_2}
\ar[u]
&
{\Gamma}
\ar@<0.3ex>[r]^-{\delta_0}
\ar@<-0.7ex>[r]_-{\delta_1} 
\ar[u]
&
{\Gamma} \ar[r]^F
\ar[u]
&
\Lambda
}
$$
which is by lemma \ref{lemma2}(ii) in each level an
strong equivalence. Using corollary \ref{intern},
we get the following equivalence of bicategories
\begin{eqnarray*}
\holim \left(
\xymatrix{
{\X\left(\Gamma_\bullet\right)}
\ar@<0.3ex>[r]^-{\delta_0^*}
\ar@<-0.7ex>[r]_-{\delta_1^*} 
&
{\X\left(\Gamma^{[2]}_\bullet\right)}
\ar@<0.9ex>[r]^{\delta_0^*}
\ar@<-0.1ex>[r]
\ar@<-1.1ex>[r]_{\delta_2^*}
&
{\X\left(\Gamma^{[3]}_\bullet\right)}
\ar@<1.3ex>[r]^-{\delta_0^*}
\ar@<0.3ex>[r]
\ar@<-0.7ex>[r]
\ar@<-1.7ex>[r]_-{\delta_3^*}
&
{\cdots}  
}
\right) \qquad\qquad\qquad\qquad\qquad
\\
\qquad\qquad\qquad\qquad\qquad
\iso \holim 
\left(
\xymatrix{
{\X\left(\Gamma_\bullet\right)}
\ar@<0.3ex>[r]^-{\delta_0^*}
\ar@<-0.7ex>[r]_-{\delta_1^*} 
&
{\X\left(\Gamma_\bullet\right)}
\ar@<0.9ex>[r]^{\delta_0^*}
\ar@<-0.1ex>[r]
\ar@<-1.1ex>[r]_{\delta_2^*}
&
{\X\left(\Gamma_\bullet\right)}
\ar@<1.3ex>[r]^-{\delta_0^*}
\ar@<0.3ex>[r]
\ar@<-0.7ex>[r]
\ar@<-1.7ex>[r]_-{\delta_3^*}
&
{\cdots}  
}
\right)
\cong \X\big(\Gamma_{\bullet}\big)
\end{eqnarray*}
The second equivalence is now a direct consequence of
Proposition \ref{bisimplicial}.
\end{Proof} 

We are now ready to take the final step and prove
theorem \ref{2.13} for $\tau$-surjective equivalences:

\begin{Proposition}\label{morita}\mbox{} \\
Let $F: \Gamma \to \Lambda$ be a $\tau$-surjective 
equivalence of Lie groupoids.
\begin{enumerate}
\item 
If $\X$ is stack, then the functor
$F^*: \X(\Lambda) \to \X(\Gamma)$ is an equivalence of
bicategories. 
\item If $\X$ is a prestack, then the functor
$F^*: \X(\Lambda) \to \X(\Gamma)$ is fully faithful.

\end{enumerate}
\end{Proposition}

\begin{Proof}
Theorem  \ref{HauptsatzEquivDesc} about equivariant descent
implies 
\begin{equation*}
\X\left(\Lambda_\bullet\right) \iso
\holim \left(
\xymatrix{
{\X\left(\Gamma_\bullet\right)}
\ar@<0.3ex>[r]^-{\delta_0^*}
\ar@<-0.7ex>[r]_-{\delta_1^*} 
&
{\X\left(\Gamma^{[2]}_\bullet\right)}
\ar@<0.9ex>[r]^{\delta_0^*}
\ar@<-0.1ex>[r]
\ar@<-1.1ex>[r]_{\delta_2^*}
&
{\X\left(\Gamma^{[3]}_\bullet\right)}
\ar@<1.3ex>[r]^-{\delta_0^*}
\ar@<0.3ex>[r]
\ar@<-0.7ex>[r]
\ar@<-1.7ex>[r]_-{\delta_3^*}
&
{\cdots}  
}
\right)
\end{equation*}
The preceding proposition  \ref{hauptsatz3} implies that this bicategory
is equivalent to $\X(\Gamma_\bullet)$, which shows part (i). 
The second statement is proven by a similar argument, using
part (ii) of theorem \ref{HauptsatzEquivDesc}.
\end{Proof}

\section{Proof of  theorem \ref{stackifizierung}}\label{sec:plusproof}

The central ingredient in the proof of theorem \ref{stackifizierung} 
is an explicit description of descent objects
$$ \Desc_{\X^+}(Y \uto M)= \X^+\big(\Ce(Y)\big) \,\, . $$
Instead of specializing to the \v Cech groupoid $\Ce(Y)$, 
we rather describe
$\X^+(\Gamma)$ for a general groupoid $\Gamma$. The plus construction involves
the choice of a cover of $\Gamma_0$ and a descent object for that cover.
For a cover $Y \uto \Gamma_0$, we consider the  
\textit{covering groupoid} $\Gamma^Y$ which is defined by
$$ \Gamma^Y_0 := Y \qquad \qquad\text{ and } \qquad\qquad \Gamma^Y_1 := Y \times_{\Gamma_0} \Gamma_1 \times_{\Gamma_0} Y \,\,\, .$$
By definition, the diagram
$$\xymatrix{
\Gamma^Y_1 \ar[rr]\ar[d]_{(s,t)} && \Gamma_1 \ar[d]^{(s,t)} \\
\Gamma^Y_0 \times \Gamma^Y_0  \ar[rr]^{\pi \times \pi} && \Gamma_0 \times \Gamma_0\\
}$$
is a pullback diagram; hence the map $\Pi: \Gamma_Y \to \Gamma$ is fully faithful
and thus a $\tau$-weak equivalence. All other structure on
$\Gamma^Y$ is induced from the groupoid structure on $\Gamma$. We thus have:

\begin{Proposition}\label{MoritaEquiv}\mbox{} \\
Let $\X$ be a prestack and $\Gamma$ be a groupoid. Then the bicategory
$\X^+(\Gamma)$ is equivalent to the following bicategory:
\begin{itemize}
\item Objects are pairs, consisting of a covering $Y \uto \Gamma_0$ and
an object $\calg$ in $\X\left(\Gamma^Y\right)$.

\item Morphisms between $(Y,\calg)$ and $(Y',\calg')$ consist of a common 
refinement $Z \uto \Gamma_0$ of $Y \uto \Gamma_0$ and $Y' \uto \Gamma_0$ 
and a morphism $A$ between the refined objects $\calg_Z$ and $\calg'_{Z}$ 
in $\X\left(\Gamma^Z\right)$.

\item 2-morphisms between one-morphisms $(Z,A)$ and $(Z',A')$ are described 
by pairs consisting of a common refinement $W \uto \Gamma_0$ of $Z$ and $Z'$ 
and a morphism of the refinements $A_W$ and $A'_{W}$ in $\X\left(\Gamma^W\right)$.

\item We identify 2-morphisms $(W,g)$ and $(W',g')$, if there exists a
common refinement $V \uto \Gamma_0$ such that the refined 2-morphisms
$g_V$ and $g'_{V}$ in $\X\left(\Gamma^V\right)$ are equal.
\end{itemize}
\end{Proposition}

\begin{Proof}
We describe explicitly an object of the bicategory $\X^+(\Gamma)$:
the first piece of data is an object in $\X^+(\Gamma_0)$. This is just
a covering $Y \uto M$ and 
\begin{itemize}
\item an object in the descent bicategory $\Desc_{\X}(Y \uto \Gamma_0)$.
\end{itemize}
The second piece of data is a morphism that relates the two pullbacks to
$\X^+(\Gamma_1)$. Such a morphism contains the coverings
$Y \times_{\Gamma_0} \Gamma_1 \uto \Gamma_1$ and 
$\Gamma_1 \times_{\Gamma_0} Y \uto \Gamma_1$ where one pullback is along
the source map and one pullback along the target map of $\Gamma$.

Proposition \ref{stabil} allows us to describe this morphism as a 
stable morphism on the canonical common refinement
$$ Y \times_{\Gamma_0} \Gamma_1 \times_{\Gamma_0} Y \uto \Gamma_1\,\,, $$
i.e.\
\begin{itemize}
\item  A morphism of pullbacks in $\Desc_{\X}(\Gamma_1^Y \uto \Gamma_1)$.
\end{itemize}

Further data and axioms can be transported to the canonical common
refinement:
\begin{itemize}
\item  A 2-morphism of pullbacks in $\Desc_{\X}(\Gamma_2^Y \uto \Gamma_2)$.
\item  A condition on the pullbacks in $\Desc_{\X}(\Gamma_3^Y \uto \Gamma_3)$.
\end{itemize}
Altogether, we have an object in 
$$
\holim \left(
\xymatrix{
{\Desc_{\X}(\Gamma_0^Y \uto \Gamma_0)}
\ar@<0.3ex>[r]^-{\partial_0^*}
\ar@<-0.7ex>[r]_-{\partial_1^*} 
&
{\Desc_{\X}(\Gamma_1^Y \uto \Gamma_1)}
\ar@<0.9ex>[r]^{\partial_0^*}
\ar@<-0.1ex>[r]
\ar@<-1.1ex>[r]_{\partial_2^*}
&
{\Desc_{\X}(\Gamma_2^Y \uto \Gamma_2)}
\ar@<1.3ex>[r]^-{\partial_0^*}
\ar@<0.3ex>[r]
\ar@<-0.7ex>[r]
\ar@<-1.7ex>[r]_-{\partial_3^*}
&
{\cdots}  
}
\right)
$$
This bicategory is, according to proposition \ref{hauptsatz3} equivalent to
$\X(\Gamma^Y)$. This shows our assertion for objects; the argument for
morphisms and 2-morphisms closely parallels the argument for objects.
\end{Proof} 

\begin{Remark}\label{beschreibungen} \mbox{} \\
We comment on the relation of the three equivalent descriptions of 
$G$-equivariant objects like e.g.\ bundle gerbes to objects described in the 
literature:

\begin{enumerate}
\item Definition \ref{defEquivariant}, which has the advantage of being
conceptually simple. This definition is used for action groupoids of
finite groups in \cite{Gawedzki-Reis}. 

\item 
The definition as a $G$-equivariant descent objects, using a $G$-equivariant open covering,
cf.\ corollary \ref{EqDesc}. This definition is used in the 
construction \cite{Meinrenken} of gerbes on compact Lie groups.

\item The characterization in proposition \ref{MoritaEquiv}, which has
the advantage that invariance under $\tau$-weak equivalences is almost immediate
from the definition. Such a definition is used in 
\cite{behrend2006dsa}.
\end{enumerate}
\end{Remark}

We are now ready for the proof of theorem \ref{stackifizierung}:

\begin{Proof}
We have to show that the presheaf $\X^+$ in bicategories is a stack.
We thus consider for any cover $Z \uto M$ the bicategory 
$$\Desc_{\X^+}(Z \uto M) = \X^+ \big( \Ce(Z) \big)\,\, .$$

By proposition \ref{MoritaEquiv}, this bicategory is given by objects,
morphisms and 2-morphisms on covering groupoids $\Ce(Z)^{Y}$ for covering
$Y \uto Z$. We write out such a groupoid explicitly:
\begin{eqnarray*}
 \Ce(Z)^{Y}_0 &=& Y = \Ce(Y)_0 \\\
 \Ce(Z)^Y_1 &=& Y \times_Z (Z \times_M Z) \times_Z Y \\
 &=& Y \times_M Y = \Ce(Y)_1 \,\,\,\, .
\end{eqnarray*}
We find $\Ce(Z)^{Y} = \Ce(Y)$. Thus $\Desc_{\X^+}(Z \uto M) 
= \X^+ \big( \Ce(Z) \big)$ is equivalent to the subbicategory of objects
of $\X^+(M)$ which are defined on coverings $Y \uto Z \uto M$. This
subbicategory is obviously equivalent to the bicategory $\X^+(M)$.
\end{Proof}

\renewcommand{\theDef}{\Alph{section}.\arabic{subsection}.\arabic{Def}}

\renewenvironment{Definition}{ \refstepcounter{Def}\mbox{}\\
\noindent\sl\textbf{Definition
\Alph{section}.\arabic{subsection}.\arabic{Def}} }
{\vspace{0.3cm}}

\appendix

\section{Appendix: surface holonomy}\label{app:shol}
\subsection{Oriented surface holonomy}\label{sapp:sholor}

To prepare the discussion of holonomy for unoriented surfaces, we 
briefly review the definition of holonomy for oriented surfaces.
The holonomy of a trivial bundle gerbe $\mathcal{I}_{\omega}$ 
with $\omega\in\Omega^2(\Sigma)$ over a 
closed oriented  surface $\Sigma$ is by definition
  \begin{equation}
  \hol_{\mathcal{I}_{\omega}} := \exp \!\Big( 2\pi\ii\int_{\Sigma} \omega \Big)
  \,\in \Ue\text\,{.}
  \end{equation}
If $\mathcal{I}_{\omega}$ and $\mathcal{I}_{\omega'}$ are two trivial bundle gerbes over 
$\Sigma$ such that there exists a 1-isomorphism $\mathcal{I}_{\omega} \to I_{\omega'}$, 
i.e.\ a U(1) bundle $L$, we have 
  \begin{equation}
  \int_{\Sigma} \omega' - \int_{\Sigma}\omega = \int_{\Sigma}\mathrm{curv}(L)
  \,\in \Z\, \label{A2}
  \end{equation}
and thus  the equality 
$\hol_{\mathcal{I}_{\omega}} = \hol_{\mathcal{I}_{\omega'}}$.

More generally, consider a bundle gerbe $\mathcal{G}$ with connection over 
a smooth oriented manifold $M$, and a smooth map 
  \begin{equation}
  \Phi:\quad \Sigma \to M
  \end{equation}
defined on a closed oriented surface $\Sigma$. Since 
$H^3(\Sigma,\Z)=0$ and since gerbes are classified by this
cohomology group, the pullback gerbe
$\Phi^{*}\mathcal{G}$ is isomorphic to a trivial bundle gerbe
$\mathcal{I}_{\omega}$. Hence one can choose a 
trivialization, i.e.\ a 1-isomorphism
  \begin{equation}
  \mathcal{T}:\quad \Phi^*\mathcal{G} \iso \mathcal{I}_{\omega}
  \end{equation}
and define the holonomy of $\mathcal{G}$ around $\Phi$ by
  \begin{equation}
  \hol_{\calg}(\Phi) := \hol_{\mathcal{I}_{\omega}}\,\text{.}
  \end{equation}
This definition 
is independent of the choice of the 1-isomorphism $\mathcal{T}$: 
if 
$\mathcal{T}' : \Phi^{*}\mathcal{G} \iso \mathcal{I}_{\omega'}$ is 
another trivialization, we have a transition isomorphism 
of gerbes on $\Sigma$
  \begin{equation}
  L:= \mathcal{T}' \circ \mathcal{T}^{-1}\colon~~\mathcal{I}_\omega \iso \mathcal{I}_{\omega'} \,\, \text{.}
  \end{equation}
The independence of the holonomy on the choice of trivialization
then follows by the argument given in equation (\ref{A2}).

\subsection{Unoriented surface holonomy}\label{sapp:sholunor}

Let $M$ be a smooth manifold and $\calj$ a Jandl gerbe on $M$. 
In this appendix, we discuss the definition of a holonomy for $\calj$
around an unoriented, possibly even unorientable, closed surface $\Sigma$. 
Such a definition is
in particular needed to write down Wess-Zumino terms for two-dimensional
field theories on unoriented surfaces which arise, e.g.\ as worldsheets
in type I string theories.

We will define surface holonomy for any pair consisting of a smooth
map $\varphi: \Sigma\to M$  and an isomorphism of $\mathbb{Z}/2$-bundles
\be  \xymatrix{
\calo(\varphi^*\calj) \ar^\sim[rr]\ar[dr] & & \hat\Sigma \ar[dl] \\
&\Sigma &
} 
\label{coveriso} \ee
where we denote the orientation
bundle of $\Sigma$ by $\hat\Sigma$. This is a canonically oriented 
two-dimensional manifold \cite{berger1988differential}.
In particular, the orientation bundle introduced in definition
\ref{newdef}.1 of the pulled back gerbe 
$\varphi^*\calj$ must be isomorphic to the orientation bundle of the surface.

\medskip

Let us first check that this setting allows us to recover the
notion of holonomy from appendix \ref{sapp:sholor} if the
surface $\Sigma$ is oriented.
An orientation of $\Sigma$ is just a global section of the orientation
bundle $\hat\Sigma\to \Sigma$.  Due to the isomorphism (\ref{coveriso}), such 
a global section gives a global
section $\Sigma\to \calo(\varphi^*\calj)$, i.e.\ an orientation of the
Jandl gerbe $\varphi^*\calj$. By proposition \ref{sequence}.2 an oriented 
Jandl gerbe 
amounts to a gerbe on $\Sigma$, for which we can define a holonomy as
in appendix \ref{sapp:sholor}.
We will see that the isomorphism in (\ref{coveriso}) is the correct
weakening of the choice of an orientation of a Jandl gerbe to the case
of unoriented surfaces.\\

Our first goal is to relate
this discussion to the one in \cite{ssw}. In that paper,
a smooth manifold $N$ together with an involution $k$ was considered.
This involution was not required to act freely, hence we describe the
situation by looking at the action groupoid $N\big/\big/(\mathbb{Z}/2)$. 
Since Jandl gerbes define a stack on $\Man$ and since any stack on
$\Man$ can be extended by definition \ref{defEquivariant} to a stack on 
Lie groupoids, the definition of a Jandl gerbe on the Lie groupoid 
$N\big/\big/(\mathbb{Z}/2)$ is clear. 

We now need a few facts about $\mathbb{Z}/2$-bundles on quotients. 
For transparency, we formulate them for the action of an arbitrary Lie group 
$G$. Consider a free $G$-action on a smooth manifold $N$ such that $N/G$ is 
a smooth manifold and such that
the canonical projection $N\to N/G$ is a surjective submersion. 
(This is, e.g., the case if the action of $G$ on $M$ is proper and 
discontinuous.) It is an
important fact that then $N\to N/G$ is a smooth $G$-bundle.

If we wish to generalize this situation to the case where the action of $G$
is not free any longer, we have to replace the quotient $N/G$ by the
Lie groupoid $N//G$. This Lie groupoid can be considered for a free
action action as well, and then the Lie groupoids $N/G$ and $N//G$ are 
$\tau$-weak equivalent. By theorem \ref{2.13}, the
categories of $G$-bundles over $N/G$ and $N//G$ are equivalent. 

This raises the question  whether there is a natural $G$-bundle on the
Lie groupoid $N//G$ generalizing the $G$-bundle $N\to N/G$.
In fact, any action
Lie groupoid $N//G$ comes with a canonical $G$-bundle $\Kan_G$ over $N//G$
which we describe as in remark \ref{rem:explicit}.
As a bundle over $N$, it is the trivial bundle $N \times G$, but it carries 
a non-trivial $G$-equivariant
structure. Namely $g \in G$ acts on $N \times G$ by diagonal multiplication, i.e.
$$g \cdot (n,h)  := (gn, gh)~ \text{.}$$ 

The following lemma shows that the $G$-bundle $\Kan_G$ has the desired
property:

\begin{Lemma}\mbox{}\label{lemma:neu} \\
Consider a smooth $G$-manifold with a free $G$-action
such that $N/G$ is a smooth manifold and such that
the canonical projection $N\to N/G$ is a surjective submersion. 
Then the pullback of the $G$-bundle $N\to N/G$ to the action Lie groupoid
$N//G$ is just $\Kan_G$.
\end{Lemma}

\begin{Proof}
The proof of the lemma consists of a careful unwinding of the definitions.
The most subtle aspect concerns the $G$-bundle
over $N$ contained in the pullback: this bundle is $N\times_{N/G} N\to N$ 
which has the diagonal as a canonical section.
\end{Proof}

We are now ready to define the target space structure corresponding to 
\eqref{coveriso}.

\begin{Definition} \mbox{} \label{def:orientifold}\\
An {\em orientifold background} consists of an action groupoid
$N\big/\big/(\mathbb{Z}/2)$, a Jandl gerbe $\calj$ on $N\big/\big/(\mathbb{Z}/2)$
and an isomorphism of equivariant $\mathbb{Z}/2$-bundles
\be  \xymatrix{
\calo(\calj) \ar^\sim[rr]\ar[dr] & & \Kan_{\mathbb{Z}/2} \ar[dl] \\
& N\big/\big/(\mathbb{Z}/2) &
}
\label{coveriso2} \ee
\end{Definition}

\begin{Proposition}\label{vergleichssw}~\\
An orientifold background is the same as a gerbe with Jandl structure 
from \cite[Definition 5]{ssw}. More precisely we have an equivalence of 
bicategories between the bicategory of orientifold backgrounds over 
the Lie groupoid $N\big/\big/(\mathbb{Z}/2)$ and the bicategory of gerbes 
over the manifold $N$ with 
Jandl structure with involution $k: N \to N$ given by the action of 
$-1 \in \mathbb{Z}/2$.
\end{Proposition}

\begin{Proof}
We concentrate on how to extract a gerbe with a Jandl structure from the 
orientifold
background. Let us first express from remark \ref{rem:explicit} the data of a
Jandl gerbe on the Lie groupoid $N\big/\big/(\mathbb{Z}/2)$ in terms of
data on the manifold $N$. We have
just to keep one isomorphism $\varphi=\varphi_k$ and a single coherence
2-isomorphism, for the non-trivial element $-1\in\mathbb{Z}/2$. We thus get:
\begin{itemize}
\item
A Jandl gerbe $\calj_N$ on $N$.
\item A morphism $\varphi: k^*\calj_N\to\calj_N$ of Jandl gerbes.
\item A coherence 2-isomorphism $c$ in the diagram
$$\xymatrix{
\calj_N \ar[rr]^{k^* \varphi}\ar[rrd]& & k^* \calj_N\ar[d]^{\varphi} \\
&  &\calj_N
\ar@{=>}_c (25,-3)*{}; (20,-8)*{}
}
$$
\item
A coherence condition on the 2-isomorphism $c$.
\end{itemize}
Similarly, we extract the data in the isomorphism
$$ \calo(\calj_N)\to \Kan_{\mathbb{Z}/2} $$
of $\mathbb{Z}/2$-bundles over the Lie groupoid $N\big/\big/(\mathbb{Z}/2)$
that is the second piece
of data in an orientifold background. It consists of

\begin{enumerate}
\item[(i)] 
An isomorphism
$$ \calo(\calj_N)\iso N\times \mathbb{Z}/2 $$
of $\mathbb{Z}/2$-bundles over the smooth manifold $N$.
\item[(ii)]
A commuting diagram
$$ \xymatrix{
\calo(k^*\calj_N) \ar_{k^*s}[d] \ar^{\calo(\varphi)}[rr] &&
\calo(\calj_N) \ar^s[d] \\
N\times \mathbb{Z}/2 \ar_{\id_N\times m_{-1}}[rr]&& N\times \mathbb{Z}/2 
} $$
where $m_{-1}$ is multiplication by $-1\in \mathbb{Z}/2$.
\end{enumerate}

Now the data in part (i) are equivalent to a section of the orientation bundle
$\calo(\calj_N)$, i.e.\ an orientation of the Jandl gerbe $\calj_N$.
By proposition \ref{sequence}.2, our Jandl gerbe is thus equivalent to an
ordinary gerbe $\calg$ on $N$. Part (ii) expresses the condition that
$\varphi$ is an orientation reversing morphism of Jandl gerbes.
We summarize the data: we get
\begin{itemize}
\item
A bundle gerbe $\calg$ on $N$.
\item
The odd morphism $\varphi$ gives, in the language of \cite{ssw}, a morphism
$A: k^*\calg\to \calg^*$ of bundle gerbes.
\item
Similarly, the coherence isomorphism 
$$ c: \varphi\circ k^*\varphi \Rightarrow \id $$
is in that language a 2-isomorphism
$$ A \otimes (k^* A)^* \Rightarrow \id $$
of gerbes which is expressed in \cite{ssw} by a $\mathbb{Z}/2$-equivariant 
structure on $A$. 
\item
Finally, one gets the coherence conditions of \cite{ssw}.
\end{itemize}

We have thus recovered all data of \cite[definition 5]{ssw}.
\end{Proof}

\begin{Kor} \mbox{} \\
The bicategory of Jandl gerbes $\calj$ over $\Sigma$ together with an 
isomorphism $f: O(\calj) \iso \hat\Sigma$ is equivalent to the bicategory of 
orientifold backgrounds over $\hat\Sigma\big/\big/ (\mathbb{Z}/2)$.
\end{Kor}

\begin{Proof}
Pull back along the $\tau$-weak equivalence $\hat\Sigma\big/\big/(\mathbb{Z}/2)
\to \Sigma$ gives by theorem \ref{2.13} an equivalence of bicategories
$$ \JGrbc(\Sigma)\iso \JGrbc(\hat\Sigma\big/\big/(\mathbb{Z}/2)) \,\, .$$
Concatenating $f$ with the isomorphism $\hat \Sigma \to \Kan_{\mathbb{Z}/2}$
from lemma \ref{lemma:neu} provides the second data in the definition
\ref{def:orientifold} of an orientifold background.
\end{Proof}

The formula for the holonomy $\hol_{\calj}(f)$ of such an orientifold background over $\hat\Sigma\big/\big/ (\mathbb{Z}/2)$ is given in 
\cite{ssw} and \cite[(5.9)]{ecm} along the lines of section \ref{sapp:sholor}. 
We refrain from giving details here. We then define

\begin{Definition} \mbox{} \label{holdef:end} \\
Let $M$ be smooth manifold and $\calj$ a Jandl gerbe on $M$. Let
$\Sigma$ be an unoriented closed surface. Given a smooth map
$\varphi: \Sigma\to M$ and a morphism $f: \calo(\varphi^*\calj)\to \hat\Sigma$
of $\mathbb{Z}/2$-bundles over $\Sigma$, we define the surface holonomy to be
$$ \hol_\calj(\varphi,f):= \hol_{(\varphi^*\calj)}(f) \,\, . $$
\end{Definition}

\begin{Remarks}\mbox{} \\[-1.8em]
\begin{enumerate} 
\item 
This holonomy enters as the exponentiated Wess-Zumino term in 
a Lagrangian description of two-dimensional sigma models on unoriented surfaces
with target space $M$ which are relevant e.g.\ for type I string theories.

\item 
More generally, one considers target spaces which are Lie groupoids.
If the target is a Lie groupoid $\Gamma$, the smooth map $\varphi$
has to be replaced by a Hilsum-Skandalis morphism $\Phi:
\Sigma \to \Lambda$ which is a special span of Lie groupoids
$$ 
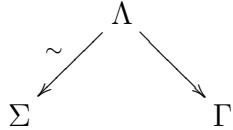
\xymatrix{
& \Lambda \ar_\sim[ld]\ar[rd] &  \\
\Sigma &&  \Gamma
}$$
where $\Lambda\to \Sigma$ is a $\tau$-weak equivalence. (For a definition
and discussion, see \cite[definition 62]{metzler2003topological}.)

Theorem \ref{2.13} ensures that the pullback along
$\Lambda\to\Gamma$ is an equivalence of bicategories. Using its inverse,
we can pull back a Jandl gerbe over $\Gamma$ along $\Phi$ to $\Sigma$.

\item 
In particular, we get in this situation a notion of holonomy 
$\hol_\calj(\Phi,f)$ for
a Hilsum-Skandalis morphism $\Phi$ and an isomorphism $f$ of 
$\mathbb{Z}/2$-bundles over $\Sigma$ as before.

\item
Consider an orientifold background, $\Gamma =N// (\mathbb{Z}/2)$. Then
each $\mathbb{Z}/2$-equivariant map $\tilde\varphi: \hat\Sigma\to N$
provides a special Hilsum-Skandalis morphism
$$ 
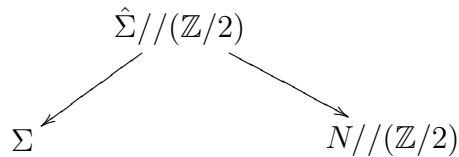
\xymatrix{
& \hat \Sigma // (\mathbb{Z}/2) \ar[ld]\ar[rd] &  \\
\Sigma &&  N // (\mathbb{Z}/2)
}$$

The pullback of $\Kan_{\mathbb{Z}/2}$ on $N// (\mathbb{Z}/2)$
to $\hat \Sigma // (\mathbb{Z}/2)$ gives again the canonical bundle which
by Lemma \ref{lemma:neu} is mapped to the $\mathbb{Z}/2$-bundle $\hat\Sigma\to\Sigma$.
Thus pulling back the isomorphism of $\mathbb{Z}/2$-bundles in the
orientifold background to an isomorphism of bundles on $\Sigma$ gives
us just the data needed in definition \ref{holdef:end} to define holonomy.

This way, we obtain holonomies $\hol_{\calj}(\tilde\varphi)$ which
have been introduced in \cite{ssw} and enter e.g.\ in orientifolds 
of the WZW models, see \cite{gsw}.

\end{enumerate}
\end{Remarks}

\bibliographystyle{alpha}
\bibliography{tncs}{}

\end{document}